\documentclass[10pt]{amsart}
\usepackage{amsthm}

\usepackage{latexsym}
\usepackage[mathscr]{eucal}
\usepackage{amssymb}
\usepackage{amsmath}
\usepackage{amsfonts}
\usepackage{a4wide}
\usepackage[mathscr]{eucal}
\usepackage[sans]{dsfont}

\DeclareMathOperator{\Hom}{Hom}
\DeclareMathOperator{\End}{End}

\DeclareMathOperator{\bimod}{bimod}
\DeclareMathOperator{\Ker}{Ker}
\DeclareMathOperator{\Aut}{Aut}
\DeclareMathOperator{\Specm}{Specm}
\DeclareMathOperator{\imm}{Im}
\renewcommand\Im{\imm}
\DeclareMathOperator{\id}{id}
\DeclareMathOperator{\St}{St}
\DeclareMathOperator\gro{growth}
\DeclareMathOperator\gr{gr}
\DeclareMathOperator\supp{supp}
\DeclareMathOperator\Sp{Spec}
\DeclareMathOperator\gkdim{GKdim}

\newcommand\ve{\varepsilon}
\newcommand\vi{\varphi}
\newcommand\la{\lambda}
\newcommand\ga{\gamma}
\newcommand{\gl}{\mathfrak{gl}}
\newcommand\Ga{\Gamma}
\newcommand\La{\Lambda}
\def\sm{\mathsf m}
\renewcommand{\k}{\Bbbk}

\newcommand\bC{\mathbb {C}}
\newcommand\bD{\mathbb {D}}
\newcommand\bN{\mathbb {N}}
\newcommand\bQ{\mathbb {Q}}
\newcommand\bZ{\mathbb {Z}}

\newcommand\cB{\mathscr{B}}
\newcommand\cG{\mathscr{G}}
\newcommand\cH{\mathscr{H}}
\newcommand\cI {\mathscr {I}}
\newcommand\cJ{\mathscr {J}}
\newcommand\cK{\mathscr{K}}
\newcommand\cL{\mathscr {L}}
\newcommand\cM{\mathscr {M}}
\newcommand\cO{\mathscr{O}}
\newcommand\cU{\mathscr {U}}

\newcommand\ds{\displaystyle}
\newcommand\myto{\longrightarrow}
\newcommand{\nad}[2]{\stackrel{#2}{#1}}
\newcommand\lr[1]{[\,#1\,]}

\newtheorem{condition}{Condition}
\newtheorem{definition}{Definition}
\newtheorem{theorem}{Theorem}[section]
\newtheorem{corollary}{Corollary}[section]
\newtheorem{lemma}{Lemma}[section]
\newtheorem{proposition}{Proposition}[section]
\newtheorem{remark}{Remark}[section]
\newtheorem{example}{Example}[section]
\newtheorem{step}{Step}

\newcommand\KVK{{_KV_K}}
\newcommand\KVL{{_KV_L}}
\newcommand\LVK{{_LV_K}}
\newcommand\LVL{{_LV_L}}
\newcommand\LWL{{_LW_L}}

\newcommand\su[2]{\overset{#2}{\underset{#1}\sum}\,}
\newcommand\osu[2]{\overset{#2}{\underset{#1}\oplus}\,}
\newcommand\bosu[2]{\overset{#2}{\underset{#1}\bigoplus}\,}
\newcommand\pru[2]{\overset{#2}{\underset{#1}\prod}\,}

\begin{document}

\author{Vyacheslav Futorny }
\author{Serge Ovsienko\  }
\title[Galois orders
\  ]{Galois orders  }
\address{Instituto de Matem\'atica e Estat\'\i stica,
Universidade de S\~ao Paulo, Caixa
Postal 66281, S\~ao Paulo, CEP
05315-970, Brasil}
\email{futorny@ime.usp.br}
\address{
Faculty of Mechanics and
Mathematics, Kiev Taras Shevchenko
University, Vla\-di\-mir\-skaya 64,
01601, Kiev, Ukraine}
\email{ovsyenko@mechmat.univ.kiev.ua }


\begin{abstract}
We introduce a new class of noncommutative rings - Galois orders,
realized as certain subrings of invariants in skew semigroup rings, and develop their
structure theory.
The class of Galois orders generalizes
classical orders in noncommutative rings and contains many classical objects, such as
the Generalized Weyl algebras, the universal enveloping algebra of
the general linear  Lie  algebra, associated Yangians and finite $W$-algebras  and certain  rings of invariant
differential operators on algebraic varieties.
  \end{abstract}

\maketitle

\tableofcontents

\section{Introduction}
\label{section-introduction} Let $\Ga$ be an integral domain and
$U\supset \Ga$ an associative noncommutative algebra over a base field $\k$. A
motivation for the study of  pairs "algebra-subalgebra" comes
 from the representation theory of Lie algebras. In particular, in the theory of Harish-Chandra modules
$U$ is the universal enveloping algebra of a reductive finite
dimensional Lie algebra $L$ and $\Ga$ is the universal enveloping  algebra of some reductive  Lie subalgebra $L'\subset L$. For instance,  the case when $\Ga$ is the universal enveloping
algebra of a Cartan subalgebra leads to the theory of
Harish-Chandra modules with respect to this Cartan algebra
- weight modules.
Another important example is a pair $(U,
\Ga)$, where $U$ is the universal enveloping algebra and $\Ga$ is  a certain maximal commutative subalgebra of $U$, called
\emph{Gelfand-Tsetlin subalgebra}. In the case $U=U(\gl_n)$ the analogs of Harish-Chandra modules - Gelfand-Tsetlin modules - were studied in \cite{dfo:gz}.
Similarly, Okounkov and Vershik (\cite{ov}) showed that representation theory of
the symmetric group $S_n$ is associated with a pair $(U, \Ga)$, where $U$ is the group algebra of $S_n$ and $\Ga$ is the
maximal commutative subalgebra generated by the Jucys-Murphy
elements.

An attempt to
understand the phenomena related to the Gelfand-Tsetlin formulae (\cite{g-ts}) was
the paper \cite{dfo:hc} where the notion of Harish-Chandra
subalgebra of an associative algebra and the corresponding notion
of a Harish-Chandra module were introduced. In particular, in \cite{dfo:hc} the categories of Harish-Chandra modules were described as categories of modules over some explicitly constructed  categories. This construction is a broad generalization of the presentation of finite dimensional associative algebras by quivers and relations. This techniques was applied to the study of  Gelfand-Tsetlin modules for $\gl_n$.

Current paper can be viewed on one hand as a development of the
ideas of \cite{dfo:hc} in the "semi-commutative case", i.e. noncommutative algebra and commutative subalgebra
 and, on the
other hand, as an attempt to understand the role of skew group
algebras in the representation theory of infinite dimensional
algebras (e.g. \cite{bl}, \cite{ba}, \cite{bavo}, \cite{ex}).
Recall, that the algebras $A_1$, $U(\mathfrak{sl}_2)$ and their quantum analogues are unified by the
notion of a \emph{generalized Weyl algebra}. Their irreducible modules are completely
described modulo classification of irreducible elements in a skew polynomial ring in one variable over a skew field.
The main property of a generalized Weyl algebra $U$ is the existence of a commutative subalgebra
$\Ga\subset U$ such that the localization of $U$ by $S=\Ga\setminus \{0\}$ is the skew polynomial algebra.
 On the other hand this technique can not be applied in case of
 more complicated algebras such as the universal enveloping algebras of simple Lie algebras of  rank $\geq 2$.

We make an important  observation that
the
Gelfand-Tsetlin formulae for $\gl_n$ define an embedding of the
corresponding universal enveloping algebra into a skew group algebra
of a free abelian group over some field of rational functions $L$
(see also \cite{khm}). A remarkable fact is that this
field $L$ is a Galois extension of the field of fractions of the
corresponding Gelfand-Tsetlin subalgebra of the universal enveloping
algebra. This fact  leads  to a
concept of  \emph{Galois orders} defined as certain subrings of
invariants in skew semigroup rings.

 We propose a notion of a "noncommutative order" as a pair $(U, \Ga)$ where $U$ is a ring,  $\Ga\subset U$ a commutative subring such that the set $S=\Ga\setminus\{0\}$ is left and right Ore subset in $U$ and the corresponding ring of fractions $\cU$ is a simple algebra (in general, $\Ga$ is not  central in $U$).  Galois orders
introduced in the paper are examples of such noncommutative orders.

Let $\Ga$ be a commutative finitely generated domain, $K$  the field of fractions of $\Ga$, $K\subset L$ a finite
Galois extension, $G=G(L/K)$  the corresponding Galois group,
$\cM\subset \Aut L$ a submonoid. Assume that $G$ belongs to the normalizer of $\cM$ in $\Aut L$ and for $m_{1},m_{2}\in \cM$ their double $G$-cosets coincide if and only if $m_{1}=gm_{2}g^{-1}$ for some $g\in G$. If $\cM$ is a group the last condition can be rewritten as
$\cM\cap G=\{e\}$. If $G$ acts on $\cM$ by conjugation then $G$ acts on the
skew group algebra $L*\cM$ by authomorphisms: $g\cdot (a m)=(g\cdot
a)(g\cdot m)$. Let $\cK=(L*\cM)^{G}$ be the subalgebra of $G$-invariants
in $L*\cM$.

We will say that an associative ring $U$ is a \emph{$\Ga$-ring}, provided  there is a fixed embedding  $i:\Ga\to U$.

\begin{definition}
A finitely generated (over $\Ga$) $\Ga$-subring $U\subset \cK$
is called a \emph{Galois $\Ga$-ring} (or \emph{Galois ring over
$\Ga$}) if $KU=UK=\cK$.
\end{definition}

 If $\Ga$ is fixed then we simply  say that $U$ is a \emph{Galois ring}.

In this case   $U\cap K$  is a maximal commutative subring in $U$ and
the center of $U$ coincides with $\cM$-invariants in  $U\cap K$
(Theorem~\ref{theorem-shat-algebras}). Moreover, the set
$S=\Ga\setminus \{0\}$ is an Ore multiplicative set (both from the
left and from the right) and the corresponding localizations
$U[S^{-1}]$ and $[S^{-1}]U$ are canonically isomorphic to
$\cK$ (Proposition \ref{proposition-Ore-condition-for-Gamma}).
Note that the
algebra $\cK$  has a canonical decomposition into a sum of
pairwise non-isomorphic finite dimensional left and right
$K$-bimodules (cf. \eqref{equation-decomp-LM-G-in-invariant}). We introduce a
class of Galois orders with some integrality conditions -
\emph{integral Galois rings}{} or \emph{Galois orders}.
These rings
satisfy some local finiteness condition and they are defined as
follows.

\begin{definition}
\label{definition-of-integral} A Galois ring $U$  over   $\Ga$ is
called \emph{right (respectively left) integral Galois ring, or Galois order over}{} $\Gamma$, if
for any finite dimensional right (respectively left) $K$-subspace
$W\subset U[S^{-1}]$ (respectively $W\subset [S^{-1}]U$), $W\cap U$
is a finitely generated right (respectively left) $\Ga$-module. A
Galois ring is \emph{Galois order}{} if it is both right
and left  Galois order.
\end{definition}
A concept of a  Galois order over $\Ga$ is a natural
noncommutative generalization of a classical notion of $\Ga$-order
in skew group ring $\cK$ since we do not require the
centrality of $\Ga$ in $U$ (cf. \cite{mcr:nnr}, chapter 5, 3.5).
 We  note the difference of our
definition from the notion of order given in \cite{mcr:nnr} (chapter
3, 1.2), \cite{hgk} (section 9).

How big is the class of Galois rings and orders? We note that any commutative
algebra is Galois over itself. If $\Ga\subset U\subset K\subset L$
and $U$ is finitely generated over $\Ga$, then $U$ is a Galois
$\Ga$-ring. If $\Gamma$ is noetherian then $U$ is an order if and
only if $U$ lies in the integral closure of $\Ga$ in $K$. In
Section~\ref{section-Bimodules} we study so-called \emph{balanced}{}
$\Ga$-bimodules. This approach, based on the bimodule theory, allows
to understand the structure and to construct Galois rings.  Another important tool in the study of
Galois rings is their Gelfand-Kirillov dimension, cf.
Section~\ref{section-Gelfand-Kirillov-dimension-of-Galois-algebras}.
Using the results of Section~\ref{section-Gelfand-Kirillov-dimension-of-Galois-algebras} we show in
Section~\ref{section-Examples-of-Galois-algebras} that the following
algebras are  integral Galois orders in corresponding skew group
rings:

\begin{itemize}
\item Generalized Weyl algebras over integral domains with infinite order automorphisms which include many classical algebras, such as $n$-th Weyl
    algebra $A_n$, quantum plane, $q$-deformed
    Heisenberg algebra, quantized Weyl algebras,
    Witten-Woronowicz algebra among the others
    \cite{ba}, \cite{bavo};

\item The universal enveloping algebra
$U(\gl_{n})$ over its Gelfand-Tsetlin subalgebra;

 \item Some  rings of invariant differential operators, e.g. symmetric and orthogonal differential operators on $n$-dimensional torus
 (cf. Section~\ref{subsection-rings-invariant-diff-operators});

 \item It is shown in  \cite{fmo},\cite{fmo1}
 that shifted  Yangians and finite $W$-al\-geb\-ras associated with  $\gl_n$ are Galois orders over corresponding Gelfand-Tsetlin subalgebras.

\end{itemize}

We emphasize  that the theory of Galois orders unifies the
representation theories of universal enveloping algebras and
generalized Weyl algebras. On one hand the Gelfand-Tsetlin formulae
give an embedding of $U(\gl_n)$ into a certain localization of the
Weyl algebra $A_m$ for $m=n(n+1)/2$
(cf. Remark~\ref{remark-embedding-gl-by} and \cite{khm}). On the other
hand the intrinsic reason for such unification is a similar hidden
skew group ring structure of these algebras as Galois orders. We believe
that the concept of  a Galois order will have a strong impact
on the representation theory of infinite-dimensional associative
algebras. We will discuss the representation theory of Galois
 rings in a subsequent  paper. Preliminary version of this paper
appeared in the preprint form \cite{fo-GaI-Pr}.

\section{Preliminaries}
\label{section-Preliminaries}

All fields in the paper contain the base  field $\k$, which
is algebraically closed of characteristic $0$. All algebras
in the paper are $\k$-algebras.

\subsection{Integral extensions}
 Let $A$ be an integral
domain, $K$ its field of fractions and $\tilde{A}$ the
integral closure of $A$ in $K$. Recall that the ring $A$ is
called \emph{normal}{} if $A=\tilde{A}$.
Let $A$ be a normal  noetherian ring, $K\subset L$ a finite
Galois extension, $\bar{A}$ the integral closure of $A$ in
$L$.

\begin{proposition}
\label{proposition-integrity}
\begin{itemize}
\item If $\tilde{A}$ is noetherian then $\bar{A}$ is
    finite over $\tilde{A}$.

\item If $A$ is a finitely generated $\k$-algebra then
    $\bar{A}$ is finite over $A$. In particular,
    $\tilde{A}$ is finite over $A$.
\end{itemize}
\end{proposition}

The following statement is probably well known but we
include the proof for the convenience of the reader.

\begin{proposition}
\label{prop-nonsingular-commutative-lifting} Let
$i:A\hookrightarrow B$ be an embedding of integral domains
with a regular $A$. Assume the induced morphism of
varieties $i^{*}:\Specm B\rightarrow \Specm A$ is
surjective (e.g. $A\subset B$ is an integral extension). If
$b\in B$ and $a b\in A$ for some nonzero $a\in A$ then
$b\in A$.
\end{proposition}

\begin{proof}
In this case $i$ induces an epimorphism of the $\Sp B$ onto
$\Sp A$. Fix  $\sm\in\Specm A$.  Assume $ab=a'\in A$. Since
the localization $A_{\sm}$ is a unique factorization
domain, we can assume that $a_{\sm}b=a_{\sm}'$, where $a_{\sm}, a'_{\sm}\in A_{\sm}$ are coprime. If $a_{\sm}$ is invertible in $A_{\sm}$ then $b\in
A_{\sm}$. If $a_{\sm}$ is not invertible in $A_{\sm}$ then there exists $P\in \Sp A_{\sm}$ such
that $a_{\sm}\in P$ and $a'_{\sm}\not\in P$. It  shows that $P$ does
not lift to $\Sp B_{\sm}$. Hence $b\in A_{\sm}$ for every
$\sm\in\Specm A$, which implies $b\in A$ (\cite{mat}, Theorem 4.7).
\end{proof}

\subsection{Skew (semi)group rings}
\label{subsection-Skew-group-rings-and-skew-group-categories}

\textit{If a semigroup $\cM$ acts on a set $S$, $\cM\times
S\xrightarrow{\vi} S$, then $\vi(m,s)$ will be denoted
 either by $m\cdot s$, or $ms$, or $s^{m}$. In particular
$s^{mm'}=(s^{m'})^{m}$, $m,m'\in \cM$, $s\in S$. By
$S^{\cM}$ we denote the subset of all $\cM$-invariant
elements in $S$. }

Let $R$ be a ring with a unit, $\cM$ a semigroup and
$f:\cM\myto \Aut (R)$ a homomorphism. Then $\cM$ acts
naturally on $R$ (from the left): $g\cdot r=f(g)(r)$ for $g\in \cM, r\in R$.
 The \emph{skew semigroup ring}{}\footnote{In a subsequent publication we will consider a more general case of the  \emph{crossed product}{} of $R$ and $\cM$.}, $R*\cM$, associated with the left action
of $\cM$ on $R$,  is a free left $R$-module,
$\ds\bigoplus_{m\in \cM} Rm$, with a basis $\cM$ and with
the multiplication defined as follows $$(r_{1}m_{1})\cdot
(r_{2}m_{2})= (r_{1}r_{2}^{m_{1}}) (m_{1} m_{2}),\quad
m_{1},m_{2}\in \cM,\ r_{1},r_{2}\in R.$$

 Assume that a finite group $G$ acts on $R$  by automorphisms and on $\cM$ by conjugation.
For every pair $g\in G$, $m\in \cM$ fix an element $\alpha_{g,m}\in R$ and define a map

\begin{equation}
\label{equation-definition-of-monomial-action}
G\times (R*\cM)\myto R*\cM,\ (g,rm)\longmapsto
\alpha_{g,m} r^g m^{g}, r\in R, m\in \cM, g\in G.
\end{equation}

This map   defines an
action of $G$ on $R*\cM$ by automorphisms if and only if $\alpha_{g,m}$ satisfy the following conditions:
\begin{align}
\label{align-conditions-on-monomial-action-in-general}
&\alpha_{g,m_{1}m_{2}}=
\alpha_{g,m_{1}}\alpha_{g,m_{2}}^{gm_{1}},\ g\in G,\ m_{1},m_{2}\in\cM;\\
&\alpha_{g_{1}g_{2},m}=\alpha_{g_{1},g_2(m)} \alpha_{g_{2},m}^{g_{1}},\ g_{1},g_{2}\in G, \alpha_{e,m}=1, m\in\cM.
\end{align}

In this case we say that the action of $G$ on $R*\cM$ is \emph{monomial}.

%

If $\alpha_{g,m}=1$ for all $g\in G$ and $m\in \cM$ then  $g(rm)=r^g m^g$ for all $g\in G$, $r\in R$, $m\in \cM$.
\textit{For simplicity
in this paper we will work with a trivial $\alpha$ but all the results remain valid for a general monomial action of $G$ satisfying the conditions above.
}

If $x\in R*\cM$ then we write it in the form
$$x=\sum_{m\in \cM} x_{m}m, $$ where only
finitely many $x_{m}\in \cM$ are nonzero. We call the
finite set
$$\supp x=\{m\in\cM|x_{m}\ne0\}$$ the
\emph{support}{} of $x$. Hence $x\in (R*\cM)^{G}$ if and only
if $x_{m^{g}}=x^{g}_{m}$ for all $m\in\cM,g\in G$.  If $x\in\
(R*\cM)^{G}$ then $\supp x$ is a finite $G$-invariant subset
in $\cM$. For $\vi\in\Aut R$ set
\begin{equation}
\label{equation-def-of-h-vi}
H_{\vi}=\{h\in G|\vi^{h}=\vi\}, \cO_{\vi}=\{\vi^{g}\mid g\in G\}.
\end{equation}

 If $G$ is a finite group
and $H$ is its subgroup then the notation $F=\ds\sum_{g\in
G/H} F(g)$ means that $g$ \textit{runs over a set of
representatives} of the quotient $G/H$ and $F(g)$ \textit{does not
depend on the choice of these representatives.} In
particular,  $F$ is well defined.

Using this agreement we denote
\begin{equation}
\label{equation-definition-a-phi}
[a\vi]:=\sum_{g\in
G/H_{\vi}}a^{g}\vi^{g}\in
(R*\cM)^{G},\quad \vi\in\cM, \ a\in R^{H_{\vi}},
\end{equation}
and set
$(R*\cM)^{G}_{\vi}=\big\{[a\vi]\,|\,a\in
R^{H_{\vi}}\,\big\}$. Then we have the following
decomposition of $(R*\cM)^G$ into a direct sum of left
(right) $R^{G}$-subbimodules
\begin{equation}
\label{equation-decomp-LM-G-in-invariant}
\begin{split}
&(R*\cM)^{G}=\bigoplus_{\vi\in
G\backslash\cM} (R*\cM)^{G}_{\vi},
\end{split}
\end{equation}
where $R^{G}$ acts on $(R*\cM)^{G}_{\vi}$ as follows
\begin{align} \label{space-structure} \gamma \cdot [a
\vi]=[(a\gamma) \vi],\ [a \vi]\cdot
\gamma= [(a\gamma^{\vi}) \vi], \
\ga\in R^{G}.
\end{align}
Thus every $x\in (R*\cM)^{G}$ can be uniquely written in the form
$\ds\sum_{\vi\in\cM\backslash G}[x_{\vi}\vi], x_{\vi}\in
R^{H_{\vi}}$.
We have for $\ga\in \Ga$
\begin{equation} \label{equation-product-of-invariants}
\begin{split}
\nonumber [a_1\vi_1]\ga[a_2\vi_2]= \sum_{\tau\in
\cO_{1}\cO_{2}}\bigg[\bigg(\sum_{\substack{g_1\in
G/H_{\vi_1}, g_2\in G/H_{\vi_2},\\
\vi_1^{g_1}\vi_2^{g_2}=\tau}}
a_1^{g_1}\ga^{\vi_{1}}a_2^{\vi_1^{g_1}g_2}\bigg)\tau\bigg].
\end{split}
\end{equation}

 For $a,b\in R^{H_{\vi}}$ denote
\begin{align}
\label{align-def-a-vi-b} &[a\vi b]=\sum_{g\in
G/H_{\vi}}a^{g}\vi^{g}b^{g}, \text{ so for } \ga\in R^{G}\text{
holds } \\ \nonumber & \ga [a\vi b]=[a\vi (b\ga^{\vi^{-1}})],\
[a\vi b]\ga=[(\ga^{\vi} a)\vi b],
\end{align}
since ${\vi}(R^{G}),{\vi^{-1}}(R^{G})\subset R^{H_{\vi}}$.
Note that in obvious way $[a\vi]=[\vi a^{\vi^{-1}}],$ $a\in R, \vi\in \Aut
R$.

\subsection{Separation actions}
\label{subsection-Separation-actions}
If $R=L$ is a field, $K\subset L$ be a finite Galois extension of fields,
$G=G(L/K)$ the Galois group and $\imath$  the canonical
embedding $K\hookrightarrow L$.
Then  $K=L^{G}$ and
\begin{equation}
\label{equation-dim-K} \dim^r_K\cK_{\vi}=
\dim^l_K\cK_{\vi}=
[L^{H_{\vi}}:K]=|G:H_{\vi}|=|\cO_{\vi}|,
\end{equation}
where $\dim^r_K,$ $\dim^l_K$ are right and left $K$-dimensions.

\begin{definition}
\label{definition-of-separating-action}
\begin{enumerate}
\item A monoid $\cM\subset \Aut L$ is called
    \emph{separating} (with respect to $K$) if
    for any $m_{1},m_{2}\in \cM$ the equality
    $m_{1}|_{K}=m_{2}|_{K}$ implies
    $m_{1}=m_{2}$.
\item An automorphism $\vi:L\myto L$ is called
    \emph{separating} (with respect to
    $K$) if the monoid generated by $\{\vi^{g}\,|\,g\in
    G\}$ in $\Aut L$ is separating.
\end{enumerate}
\end{definition}

\begin{lemma}\label{lemma-separ-monoid} Let  monoid $\cM$ be separating with respect to
    $K$. Then
\begin{enumerate}
\item \label{enum-def-sep-action-nonunit} $\cM\cap G=\{e\}.$

\item \label{enum-def-of-sep-action-act} For any
    $m\in\cM,m\ne e$ there exists $\ga\in K$ such that
    $\ga^{m}\ne\ga$.

\item\label{enum-def-of-sep-action-galois} If $G
    m_{1}G=G m_{2} G$ for some $m_{1},m_{2}\in\cM$,
    then there exists $g\in G$ such that
    $m_{1}=m_{2}^{g}$.

\item\label{enum-def-of-sep-action-group-case} If $\cM$ is a group, then the statements \eqref{enum-def-sep-action-nonunit}, \eqref{enum-def-of-sep-action-act}, \eqref{enum-def-of-sep-action-galois} are equivalent and each of them imply that $\cM$ is separating.
\end{enumerate}
\end{lemma}

\begin{proof}
We prove the statement \eqref{enum-def-of-sep-action-galois}, other statements are trivial.  $G  m_{1}G=G m_{2} G$ if and only if for some $g,g'\in G$ holds $m_{1}^{g}=m_{2}g'.$ Then  $m_{1}^{g}$ and $m_{2}$ acts in the same way on $K$, hence $m_{1}^{g}=m_{2}.$
\end{proof}

Let $\jmath:K\hookrightarrow L$ be an embedding. Denote
$\St(\jmath)=\{g\in G|g\jmath=\jmath\}$.

\begin{lemma}
\label{lemma-separating-action-and-stabilizer} Let
$\vi\in\cM$, $\jmath=\vi\imath$.  Then
\begin{enumerate}
\item\label{enumerate-H-vi-equals-stab-vi}
If $\vi$ is  separating, then  $H_{\vi}=\St(\jmath)$.
\item\label{enumerate-H-vi-stabilize-K-vi-K}
    $K\vi(K)=L^{\St(\jmath)}$, in particular, if $\vi$ is separating $K\vi(K)=L^{H_{\vi}}$.
\end{enumerate}
\end{lemma}

\begin{proof}
Obviously $H_{\vi}\subset \St(\jmath)$. Conversely, if
$g\vi\imath=\vi\imath$, then $\vi^{-1}g\vi\imath=\imath$,
hence $\vi^{-1}g\vi=g_{1}\in G$ and $\vi^{-1}(g\vi
g^{-1})=g_{1}g^{-1}$. Thus $\vi$  and $g\vi g^{-1}$
coincide on $K$, implying $g\vi g^{-1}=\vi$ and
\eqref{enumerate-H-vi-equals-stab-vi}.
Note that $g\in G(L/K\vi(K))\cap G$ if and only if $g|_{\vi(K)}=\id$ (i.e. $g\in \St (\jmath)$), implying
\eqref{enumerate-H-vi-stabilize-K-vi-K}.
\end{proof}


\section{Bimodules}
\label{section-Bimodules}

\subsection{Balanced bimodules}
\label{subsection-Bimodules-and-homomorphisms}

For commutative $\k$-algebras $A$ and $B$ we will denote by
$(A-B)-\bimod$ the category of finitely generated
$A-B$-bimodules. If $A=B$ we will simply write $A-\bimod$.

\begin{proposition}
\label{lemma-krull-schmidt-for-bimodules} Let $K\subset L$
be a finite field extension. The full subcategories of $K-\bimod$,
$(K-L)-\bimod$ or $(L-K)-\bimod$ consisting of objects,
which are finite dimensional as left or as right modules
are
 Jordan-Hoelder and Krull-Schmidt categories.
\end{proposition}

\begin{proof} It follows from the finiteness of the length of the objects of these categories.
\end{proof}

\textit{In this section all bimodules over  fields are assumed to be finite
dimensional from both sides and $\k$-central (unless the
contrary is stated).}

\begin{definition}
\label{definition-bimodule-induced-by-homomorphism} A
homomorphism of  algebras $\vi:A\to B$  endows $B$ with the  structure of $B-A$-bimodule
$B_{\vi}$  such that for $a\in A, b\in B$, $b'\in B_{\vi}$ holds  $b\cdot
b' \cdot a= bb'\vi(a).$
\end{definition}

\begin{remark}
\label{remark-bimod-produces-hom-and-tensor-corresponds-composition}
\begin{enumerate}
\item\label{enumerate-bimodule-produces-hom}

In opposite,  an $B-A$ bimodule
$V$, which is free of rank $1$ from the left, defines a
homomorphism $\vi=\vi_{V}:A\to B$ by $va=\vi(a)v $, where
$v\in V$ is a right free generator of $V$.

\item\label{enumerate-tensor}
If $\vi:A\myto B$ and $\psi:B\myto C$ are homomorphisms of algebras
then there exists an isomorphism of $C-A$-bimodules $$C_{\psi} \otimes_{B} B_{\vi}\simeq  C_{\psi\vi},\ c\otimes b \longmapsto c \psi (b), c\in C,b\in B.$$
\end{enumerate}
\end{remark}

Let $K\subset L$ be an extension and $\imath_{K}$  the canonical embedding
$K\subset L$. We will write  $\imath$  instead of  $\imath_{K}$ when the field $K$ is
fixed. If
$V=\KVK$ is a $K$-bimodule then  denote $\KVL= V\otimes_K L
$, $\LVK=L\otimes_K V$ and $\LVL=L\otimes_K \KVL$.

Let $K\subset L$ is a Galois extension with the Galois group
$G=G(L/K)$, then
the group $G\times G$ acts on $\LVL$  as
$$(g_{1},g_{2})\cdot (l_{1}\otimes v\otimes l_{2})\longmapsto l_{1}^{g_{1}}\otimes v\otimes l_{2}^{g_{2}^{-1}}, (g_{1},g_{2})\in G\times G, v\in V, l_{1},l_{2}\in L$$ by automorphism of $K$-bimodules. The $K$-bimodule of invariants is canonically isomorphic to $V$.  If we restrict the action of $G\times G$ to the action of $G$ from the left (from the right), by automorphisms of $K-L$ ($L-K$)  bimodules, then the invariants will be $\KVL$ ($\LVK$).

Analogously, $G$ acts naturally from the left on the $L-K$-bimodule
$\LVK$ by automorphisms of $K$-bimodule, $$g\cdot (l\otimes v)\longmapsto l^{g}\otimes v, g\in G, v\in V, l\in L \text{ and } (\LVK)^{G}\simeq \KVK.$$

Assume that the right action of $K$ on $V$ is
$L$-diagonalizable from the left. It means $\LVK$ splits into
a sum of $L-K$-bimodules, which are one dimensional as
right $L$-modules. By Remark \ref{remark-bimod-produces-hom-and-tensor-corresponds-composition}, \eqref{enumerate-bimodule-produces-hom}  such one dimensional $L-K$-bimodule is of the form $L_{\jmath}$ for some field
embedding $\jmath:K\to L$.

\begin{definition}
\label{definition-balanced-bimodule} A $K$-bimodule $\KVK$
is called \emph{$L$ -{ balanced}}{} over a finite Galois
extension $K\subset L$ if $\LVL$ is a direct sum of
one-dimensional from the left and from the right $L$-bimodules, i.e. bimodules of the form $L_{\vi}$ for $\vi\in \Aut L.$ A $K$-bimodule $\KVK$ is called \emph{balanced}{} if it is $L$-\emph{balanced} over
some finite Galois extension $K\subset L$.
\end{definition}

\subsection{Monoidal category of balanced bimodules}
\label{subsection-Tensor-category-of-balanced-bimodules}
Denote by $K-\bimod_{L}$ the full subcategory in $K-\bimod$
consisting of all $L$-balanced $K$-bimodules.

\begin{remark}
\label{remark-simple-L-L-balanced}
The category $L-\bimod_{L}$ is by definition semisimple and its isoclasses of simples are represented by the bimodules $L_{\vi}$, were $\vi:L\myto L$ is an automorphism.
\end{remark}

\begin{theorem}
\label{theorem-tensor-category-of-balanced}
 The category $K-\bimod_{L}$ is
     an abelian semisimple monoidal category.

\end{theorem}

\begin{proof}
Note that by Remarks \ref{remark-bimod-produces-hom-and-tensor-corresponds-composition}, \eqref{enumerate-tensor} and by  Remark \ref{remark-simple-L-L-balanced} above the category $L-\bimod_{L}$ satisfies the theorem.

Let $V,W$ be $L$-balanced $K$-bimodules, $\pi: V\myto W$ an $K$-bimodule
epimorphism, $\pi_L:\LVL\myto {}_{L}W{}_{L} $  the induced
epimorphism of $L$-bimodules.
Since $G$ acts trivially on $K$ the map $\pi_{L}$ is a homomorphism of $ (K\otimes_{\k}K)[G\times G]$-bimodules.

On the other hand $p_L$ admits the right
inverse $L-L$-bimodule monomorphism
$$s_{L}:{}_{L}W{}_{L}\myto \LVL,\ p_{L} s_{L}=\id_{{}_{L}W{}_{L} }.$$
Since $G$ acts trivially on $K$ for every $g=(g_{1},g_{2})\in G\times G$ the
morphism $$g s_{L} g^{-1}: {}_{L}W{}_{L}\myto \LVL,\quad l_{1}\otimes w \otimes l_{2}\
\longmapsto g_{1}\cdot s_{L}(l_{1}^{\,g_{1}^{-1}}\otimes w \otimes l_{2}^{g_{2}})\cdot g_{2}^{-1}$$ are
$K$-bimodule homomorphisms.  Then the $K$-bimodule homomorphism
$$\sigma_{L}=\dfrac{1}{|G|^2}\sum_{g\in G\times G} g s_{L}g^{-1}$$
commutes with the action $G\times G$, hence both $\sigma_{L}$ and $\pi_{L}$ are $(K\otimes_{\k}K)[G\times G]$-bimodule homomorphisms. We have $$p_{L} \sigma_{L}=\ds\dfrac{1}{|G|^2}\sum_{g\in
G\times G} p_{L}g s_{L}g^{-1}=\ds\dfrac{1}{|G|^2}\sum_{g\in
G\times G} gp_{L} s_{L}g^{-1}=\id_{{}_{L}W{}_{L} }.$$ Since
$\sigma_{L}$  maps ${}_{L}W{}_{L}^{G\times G}$ to
$\LVL^{G\times G}$, it induces a $K$-bimodule
homomorphism $\sigma: W\myto V$, which splits $p$. Hence  $K-\bimod_{L}$ is semisimple.

Consider the standard $K$-bimodule monomorphism $$i:V\otimes_{K}W\myto
V\otimes_{K}L\otimes_{K}W,\quad v\otimes w\longmapsto v\otimes 1 \otimes w.$$
Then the induced $L$-bimodule
homomorphism
$$_{L}(V\otimes_{K}W)_{L}\myto
L\otimes_{K}V\otimes_{K} {L}\otimes_{K}W\otimes_{K}L\simeq \LVL\otimes_{L}\LWL,$$  is a monomorphism.
Since $_{L}{V}_{L}$ and $_{L}{W}_{L}$ are
isomorphic to the sums of simple one-dimensional
$L$-bimodules, the same is true for their tensor product
over $L$ and for its subbimodule $_{L}(V\otimes_{K}W)_{L}$.

Note also that $K$ is a weak unit with respect to $\otimes_{K}$ in $K-\bimod_{L}$.
\end{proof}

\subsection{Simple balanced bimodules}
\label{subsection-Simple-balanced-bimodules} In this
section we describe simple objects in
$K-\bimod_{L}$.

\begin{lemma}
\label{lemma-on-balanced-sufficient}
Let $K\subset L$ be a Galois extension, $G=G(L/K)$.
\begin{enumerate}
    \item\label{enumerate-on-tensor-product-L-vi-over-K}(\cite{dk}, Ch. 5.1) If for a field $F$ holds $K\subset F\subset L$, $H=G(L/F)$ and $i_{F}:F\hookrightarrow L$ is the canonical embedding, then as $L-F$-bimodule
     $$ L\otimes_{K} F\simeq \bosu{g\in G/H}{} L_{g\, i_{F}}, \text{
in particular } L\otimes_{K} L \simeq \bosu{g\in G}{} L_{g}.$$

    \item\label{enumerate-on-balanced-sufficient} A $K$-bimodule $V$ is $L$-balanced if and only if the $L-K$-bimodule $\LVK$ is a direct sum of modules of the form $L_{\vi\imath}$, $\vi\in\Aut L.$

    \item\label{enumerate-left-right-dimension}
The right and the left $K$-dimensions of a balanced bimodule coincide.
\end{enumerate}
\end{lemma}

\begin{proof}
To prove the statement \eqref{enumerate-on-tensor-product-L-vi-over-K} we present $F$ as a simple extension $F=K[\alpha]$, $\alpha\in F$. Let $f(X)$ be a minimal polynomial of $\alpha$ over $K$, $\alpha=\alpha_{1},\dots,\alpha_{k}\in L$  all roots of $f(X)$.  Then $F\simeq K[X]/(f(X))$ and
$$ L\otimes_{K} F\simeq L\otimes_{K}K[X]/(f(X))\simeq  L[X]/(f(X))\simeq\prod_{i=1}^{k}L[X]/(X-\alpha_{i}).$$ The right $F$-module structure on $L[X]/(X-\alpha_{i})$ is defined by multiplication on $X$, that proves \eqref{enumerate-on-tensor-product-L-vi-over-K}.

In \eqref{enumerate-on-balanced-sufficient} we prove firstly the statement ``if''. Applying Remark \ref{remark-bimod-produces-hom-and-tensor-corresponds-composition}, \eqref{enumerate-tensor} we obtain the following isomorphisms of $L-K$-bimodules, which proves the statement.
\begin{align*}
&L\otimes_{K}L_{\vi\imath}\simeq L\otimes_{K}(L\otimes_{L} L_{\vi\imath})\simeq(L\otimes_{K}L)\otimes_{L}  L_{\vi\imath}\simeq \\
&\big(\bosu{g\in G}{} L_{g}\big)  \otimes_{L}   L_{\vi\imath}\simeq\bosu{g\in G}{} L_{g\vi \imath}.
\end{align*}

Prove the statement ``only if''. If $\LVL\simeq\bosu{\vi\in S}{}L_{\vi}^{d_{\vi}},$ $S\subset\Aut L,$ $d_{\vi}>0$ as $L-L$-bimodule, then as $L-K$-bimodule it is isomorphic to $\bosu{\vi\in S}{}L_{\vi\imath}^{d_{\vi}}$. In particular, $\LVL$  is a semisimple $L-K$-bimodule.
 Note, that $\LVK$ can be identified with $ (\LVL)^{\{e\}\times G}$ which is a  $L-K$-submodule in $\LVL$. Hence $\LVK$ as a subbimodule of the semisimple $L-K$-bimodule $\LVL$ is a direct sum of some $L_{\vi\imath},$ $\vi\in S.$

 The statement \eqref{enumerate-left-right-dimension} follows from the definition.
\end{proof}

\begin{lemma}
\label{lemma-simple-balanced-preliminary}
Let $\vi\in \Aut L$, $\jmath=\vi\imath$,  $H=\St(\jmath)$.
\begin{enumerate}
\item\label{enumerate-prelim-simple-balanced-action}
The canonical action of  $H$ on $L$ is an action by $L^{H}-K$-bimodule automorphisms on $L$-bimodule $L_{\vi}$ and on $L-K$-bimodule $L_{\jmath}.$

\item\label{enumerate-prelim-simple-balanced-v-vi}
Let $j:K\to L^{H}$ induced by $\jmath$ embedding,  $V(\vi)=L_{\vi}^{H}$. Then $V(\vi)$ as $L^{H}-K$ subbimodule in $L_{\vi}$ is isomorphic to $L_{î}$. In particular, $V(\vi)$ is a simple $L^{H}-K$-bimodule.
\end{enumerate}
\end{lemma}

\begin{proof} Let $l\in L_{\vi},$ or $l\in L_{\jmath},$ $l_{1}\in L^{H},$ $k\in K$ and $''\cdot''$ is the bimodule action on $L_{\vi}$. Then for $h\in H$ holds
\begin{align*}
&(l_{1}\cdot l\cdot k)^{h}=(l_{1}l \vi(k))^{h}=l_{1}^{h} l^{h} h\vi(k)=l_{1}l^{h} \vi(k)=l_{1}\cdot l^{h} \cdot k,
\end{align*}
which proves the statement \eqref{enumerate-prelim-simple-balanced-action}. Further,
$V(\vi)$ is $L^{H}-K$-bimodule by \eqref{enumerate-prelim-simple-balanced-action}. Other statements of \eqref{enumerate-prelim-simple-balanced-v-vi} are obvious.
\end{proof}

$V(\vi)$ has a structure of $K$-bimodule since $K\subset L^{H}$. It turned out, that $V(\vi),$ $\vi\in\Aut L$ cover all simples in $K-\bimod_{L}.$

\begin{theorem}
\label{proposition-construction-of-balanced}
\begin{enumerate}
\item
\label{enum-constructed-is-balanced}
$L\otimes_{K}V(\vi)\simeq \ds\bigoplus_{g\in G/H}
L_{g\vi\imath}$ as a $L-K$-bimodule, i.e. $V(\vi)$ is $L$-balanced.

\item \label{enum-v-vi-is-simple-balanced} $V(\vi)$ is
    a simple $K$-bimodule.

\item \label{enum-all-simple-balanced} Any simple
    object in  $K-\bimod_{L}$ is isomorphic
    to  $V(\vi)$ for some $\vi\in\Aut L$.

\item \label{enum-isomorphic-simple-balanced} Let
    $\vi,\vi'\in \Aut L$. Then $V(\vi)\simeq V(\vi')$
    if and only if $\,G\vi|_{K}=G\vi'|_{K}$,
    equivalently $G\vi G=G\vi' G$.

\item \label{enumerate-K-vi-K} Assume $\vi\in \cM$ for a separating
    monoid $\cM\subset \Aut L$,  $a\in L^{H}$,
    $v=[a\vi]\in \cK$, \eqref{equation-definition-a-phi}. Then $KvK\simeq V(\vi)$ as
    $K$-bimodule.
\end{enumerate}
\end{theorem}

\begin{proof}
Consider $V(\vi)$ as  $L^{H}-K$ bimodule.
Using Lemma \ref{lemma-on-balanced-sufficient},  \eqref{enumerate-on-tensor-product-L-vi-over-K} and Remark \ref{remark-bimod-produces-hom-and-tensor-corresponds-composition}, \eqref{enumerate-tensor} we obtain the following isomorphisms of $L-K$-bimodules, which proves \eqref{enum-constructed-is-balanced}
\begin{align*}
& L\otimes_{K}V(\vi)
= L\otimes_{K}L^{H}_{\jmath}\simeq  L\otimes_{K}( L^{H}\otimes_{L^{H}}L^H_{\jmath}) \simeq (L\otimes_{K}L^{H})\otimes_{L^{H}}    L^H_{\jmath}\simeq \\& \text{  }(\bosu{g\in G/H}{}L_{g})\otimes_{L^{H}} L^H_{\jmath}\simeq \bosu{g\in G/H}{}(L_{g}\otimes_{L^{H}}L^H_{\jmath})\simeq
 \bosu{g\in G/H}{}L_{g\jmath}.
\end{align*}
To prove the simplicity of $V(\vi)$ consider any nonzero
$x\in L^{H}$. Then $K\cdot x\cdot K=x\vi(K)K= L^{H}$,  by Lemma \ref{lemma-separating-action-and-stabilizer},\eqref{enumerate-H-vi-stabilize-K-vi-K},  implying
\eqref{enum-v-vi-is-simple-balanced}.

Now we prove \eqref{enum-all-simple-balanced}.
Let $V$ be a simple  $L$-balanced $K$-bimodule. We divide the proof in the the following steps. If $A$ is a $\k$-algebra, then in the proofs below instead of s structure of $A-K$-bimodule we will use the corresponding structure of left $A\otimes_{\k}$-module.

\begin{step}
\label{step-LstarG-structure}
The equality $( l'g\otimes k)\cdot
(l\otimes v)=l' l^g\otimes kv$, $k\in K$, $g\in G$, $l,
l'\in L$, $v\in V$ endows  $\LVK$ with the structure of a simple left
$(L*G)\otimes_{\k}K $-module.
\end{step}

The correctness of $(L*G)\otimes_{\k}K $-module structure is checked immediately. To prove the simplicity consider  $0\ne x\in \LVK$, $x=\su{g\in G}{}l_{g}\otimes v_{g}$, where $v_{g}\in V,g\in G$ and  $\{l_{g}\mid l\in L,g\in G\}$ is a normal $K$-basis of $L$. Consider $g'\in G$ such that $v_{g'}\ne0$. By the theorem of independence of characters the maps $w_{g}:G\to L$, $w_{g}(g_{1})=l_{g_{1}g},$ $g\in G$ form a basis in the $L$-vector space of maps $G\to L$. Hence
there exist $\su{g\in G}{} \la_{g}g\in L*G,$ $\la_{g}\in L$, such that $$ \big(\su{g\in G}{} \la_{g}g\big)\cdot x=\su{g\in G}{} \bigg(\su{g_{1}\in G}{}\la_{g_{1}} l_{g_{1}g} \bigg)\otimes v_{g}=1\otimes v_{g'},$$ which obviously generates $V$ as $K$-bimodule and $\LVK$ as $L-K$-bimodule.

\begin{step}
\label{step-LstarG-summsnds}
 $ \LVK\simeq \osu{g\in G/H}{} L_{g\jmath}^{d}$  for some $d\geq1$, where $\jmath=\vi\imath$ for some $\vi\in \Aut L$. Besides every $L_{g\jmath}$  is a simple $(L*H)\otimes_{\k}K$-submodule in $\LVK$, where $H=\St(g\jmath)$.
\end{step}

By definition   $\LVK\simeq  \bosu{\jmath\in S}{} L_{\jmath}^{d_{\jmath}}$ as a $L-K$-module for some pairwise non-isomorphic $L_{\jmath}$.
Since $g(L_{\jmath})\simeq L_{g\jmath}$ and $\LVK$ is simple as $(L*G)\otimes_{\k}K $-module, we have
$ \LVK\simeq \osu{g\in G/H}{} L_{g\jmath}^{d}$ as a $L-K$ bimodule.
The  $L-K$-subbimodule $L_{g\jmath}$  of $\LVK$ is $H$-invariant,  hence it is a $(L*H)\otimes_{K}K$-module, $H=H_{g\jmath}$.
Besides, $L_{g\jmath}$ is irreducible as $L-K$-bimodule.

\begin{step}
\label{step-d-equals-1}
$d=1$.
\end{step}
Note that $(L*G)\otimes_{\k}K $ is a free right $(L*H)\otimes_{\k}K$-module of rank $[G:H]$.
The canonical
embedding of $(L*H)\otimes_{\k}K$-modules
$L_{\jmath}\hookrightarrow \LVK$ induces a homomorphism of
$(L*G)\otimes_{\k}K $-modules
$$\Phi:(L*G)\otimes_{L*H} L_{\jmath}\myto \LVK,$$
which is an epimorphism, since
$\Phi\ne0$ and $\LVK$ is simple.
On the other hand for the left $K$-dimensions $\dim^l_K$ holds
$$\dim^l_K (L*G\otimes_{L*H} L_{\jmath})=[L:K][G:H],\quad\dim^l_{K}\LVK=d[L:K][G:H].$$
Hence, $d=1$ and $\Phi$ is an isomorphism.

\begin{step}
\label{step-Phi-iso}
The mapping $$\psi:K[G]\times L_j\myto (L*G)\otimes_{L*H} L_{\jmath},\quad (kg,l)\longmapsto
kg\otimes l,\ k\in K,g\in G, l\in L_{j}$$ induces an isomorphism of left $K[G]\otimes_{\k}K$-modules $$\Psi:K[G]\otimes_{K[H]} L_j\myto (L*G)\otimes_{L*H} L_{\jmath}.$$
\end{step}

Indeed,
$\psi$ is $K[H]$-bilinear and commutes with the action of $K[G]$ from the left and with the action of $K$ from the right.
Again a comparison of $K$-dimensions implies the statement.

\begin{step}
\label{step-finish}
$V\simeq V(\vi).$
\end{step}

Steps \eqref{step-d-equals-1} and \eqref{step-Phi-iso} shows, that the composition $$\Psi\circ\Phi: K[G]\otimes_{K[H]}L_{\jmath}\myto \LVK$$ is an isomorphism of $K[G]\otimes_{\k}K$-modules.
By the Frobenius reciprocity for  left $K[H]$-module $L_{\jmath}$ we obtain the chain of $K$-bimodule isomorphisms
\begin{align*}
&V\simeq (\LVK)^{G}\simeq(K[G]\otimes_{K[H]}L_{\jmath})^G\simeq \Hom_{K[G]}(K, K[G]\otimes_{K[H]}L_{\jmath})\simeq\\
&\Hom_{K[G]}(K, \Hom_{K[H]}(K[G],L_{\jmath}))\simeq\Hom_{K[H]}(K[G]\otimes_{K[G]}K, L_{\jmath})
\simeq \Hom_{K[H]}(K, L_{\jmath})\simeq
L_{\jmath}^H.
\end{align*}

It leaves to prove \eqref{enumerate-K-vi-K}.
Assume $V(\vi)\simeq V(\vi')$.
Then $L\otimes_{K}V(\vi)\simeq L\otimes_{K}V(\vi')$. Hence
from \eqref{enum-constructed-is-balanced}, $\vi'\imath =
g\vi\imath$ for some $g\in G$ and thus $G\vi \imath =G\vi'
\imath$ and $G\vi|_{K}=G\vi'|_{K}$ The converse statement easily follows.

Using  \eqref{align-def-a-vi-b}  and Lemma
\ref{lemma-separating-action-and-stabilizer},
\eqref{enumerate-H-vi-stabilize-K-vi-K} we obtain
$K[a\vi]K=[K\vi(K)a\vi]=[L^{H}a\vi]$ which immediately implies the isomorphism
$[L^{H}a\vi]\simeq V(\vi)$ and hence the last statement.
\end{proof}

\subsection{Grotendieck ring of category  balanced bimodules and Hecke algebra}
\label{subsection-Grotendieck-ring-of-category-of-balanced-bimodules}

Let $K_{0}(K,L)$ be the  Grothendieck ring
of the category $K-\bimod_{L}$ and for $V\in K-\bimod_L$ $[V]$ the class of $V$ in $K_{0}(K,L)$.
Theorem \ref{proposition-construction-of-balanced} shows  that simple $L$-balanced  $K$-bimodules can be enumerated by the double cosets $G\vi G$ or by the $G$-orbits $G \vi \imath$. We show that the ring structure on $K_{0}(K,L)$ is closely related with some Hecke algebra (Corollary \ref{corollary-rational-Grothendieck-isomorphic-Hecke}).

To calculate in $K_{0}(K,L)$ we need some preliminaries.
A \textit{family of elements} $S$ of a set $T$ is the mapping $S:\cI\myto T$, where $\cI$ in the set of indices. If the group $G$ acts on $\cI$ and $T$, then we say  $S$ is $G$-invariant provided that $S$ is a  map of $G$-sets. To simplify the notation we will write $i$ instead
of $S(i), i\in\cI$.
By $S/G$ we denote the induced map of factorsets $S/G:\cI/G\myto T/G$.
In particular, $S/G$ is a family of elements of $T/G$, indexed by $\cI/G$.

Denote $\Hom_{\k-fiels}(K,L)$ the set of all field $\k$-embeddings $K\to L$, and   $$\cB(K,L)=\{S\mid S\subset \Hom_{\k-fiels}(K,L),\ |S|<\infty, GS=S \}.$$
Then by Lemma \ref{lemma-on-balanced-sufficient}, \eqref{enumerate-on-balanced-sufficient} we can correspond to a finitely generated  balanced $K$-bimodule $V$  a $G$-invariant family $S_{V}:\cI_{V}\myto \Hom_{\k}(K,L)$, such that $\LVK\simeq \bosu{\tau\in
\cI_{V}}{}L_{S_{V}(\tau)}$. In obvious way factorization by $G$ induces the family $s_{V}=S_{V}/G: \cI_{V}/G\myto \cB(K,L).$ Obviously, the image of $s_{V}$ defines the $K$-bimodule $V$ uniquely up to isomorphism and we can write $\LVK\simeq \bosu{\tau\in
\cI_{V}/G}{}L_{s_{V}(\tau)}$.

In particular, by Theorem \ref{proposition-construction-of-balanced} \eqref{enum-constructed-is-balanced}, we can choose  $\cI_{V(\vi)}$ the set $G/H$, $S_{V}(gH)=g\vi.$
Then $\cI_{V(\vi)}/G$ is one-element and and the image of $s_{V}$ is the subset $\{g\vi\mid g\in G/H,\ H=\St(\vi \imath)\}.$
A double coset $C=G\vi G\in G\backslash \Aut L/G$ defines  an
$$b_{C}=b_{\vi}=\sum_{\psi\in C}\psi=
\sum_{g\in G/H_{\vi}} \sum_{\tau\in
g\vi G}\tau\in \bQ[\Aut L].$$ If
$x=\su{\vi\in G\backslash \Aut L/G}{} n_{\vi}b_{\vi}\in \bQ[\Aut L]$,
$n_{\vi}\in \bN$, then denote
$V(x)=\bosu{\vi\in G\backslash \Aut L/G}{}
V(\vi)^{n_{\vi}}.$ In particular,
$V(b_{\vi})\simeq V(\vi)$.

\begin{corollary}\label{corollary-multoplicities-of-simples}
Let $V$ be an object of $K-\bimod_{L}$ and in the notation above $V\simeq
\bosu{\tau\in
\cI_{V}/ G}{}V(s_{v}(\tau)).$

\begin{enumerate}
\item\label{enumerate-multiplicity-in-balanced}
For $\vi\in\Aut L$  the multiplicity $n_{\vi}$ of
$V(\vi)$ in $V$ is given by
$$n_{\vi}= \su{\tau\in
\cI_{V},\ S_{v}(\tau)=\vi\imath}{}\dfrac{|\St(\vi\imath)|}{|G|}, $$
\item\label{enumerate-class-in-Grothendieck-group}
$[V]=\su{\tau\in\cI_{V}}{} \dfrac{|\St(S_{V}(\tau))|}{|G|}[V(S_{V}(\tau))].$
\end{enumerate}
\end{corollary}

\begin{proof}
The statement \eqref{enumerate-class-in-Grothendieck-group} follows from \eqref{enumerate-multiplicity-in-balanced}
The proof follows from Theorem \ref{proposition-construction-of-balanced}, \eqref{enum-constructed-is-balanced}.
\end{proof}

 Recall,  if $G_{1}$ is a
group, $G\subset G_{1}$ is a finite subgroup and $A$ is
a commutative ring, then the Hecke algebra
$\cH_{A}(G_{1};G)\subset A[G_{1}]$ is a free module
over $A$ with a basis $h_{G\vi G}$ labeled by double cosets
in $G\backslash G_1/G$. For details on Hecke algebras
we refer to \cite{kr}. We will need the following result
from \cite{kr} (Theorem 1.6.6) slightly adapted to our
conditions.

\begin{theorem}
\label{theorem-hecke-algebras-references}
Let $\Omega=\Aut L$. Then
\begin{enumerate}
\item $e_G=\ds\dfrac{1}{|G|}\sum_{g\in G}g$ is an
    idempotent in the group algebra $\bQ[\Omega]$.

\item One has
$e_G \vi e_G=\dfrac{|H_{\vi}|}{|G|^2}b_{\vi}$ for all $\vi\in \Omega$ and
$e_G \bQ[\Omega] e_G$ becomes a subalgebra of $\bQ[\Omega]$
with $e_G$ as its identity element.

\item The mapping
$\Phi: \cH_{\bQ}(\Omega;G)\myto e_G \bQ[\Omega] e_G\subset \bQ[\Omega],$ where
$$\sum_{\vi\in G\backslash \Omega/G}n_{\vi}h_{G\vi G}\longmapsto \dfrac{1}{|G|}\sum_{\vi\in G\backslash \Omega/G}b_{\vi}$$ is an isomorphism of $\bQ$-algebras.
\end{enumerate}

\end{theorem}

We will identify the Hecke algebra $\cH_{\bQ}(\Omega;G)$
with $\Im(\Phi) \subset \bQ[\Omega]$. Given $\vi,\psi\in
\Aut L$, introduce an equivalence relation $\sim$ ($=\sim(\vi,\psi)$) on $G$ as
follows:
$$g\sim g'\textrm{ if and only if } G\vi g\psi G=G\vi g'\psi G.$$
For $\vi\in\Aut L$ denote by $H_{\vi}=\St(\vi\imath)$,
where $\imath:K\hookrightarrow L$ is the canonical
embedding.

\begin{theorem}
\label{theorem-tensor-product-of-L-balanced}
Let $\vi,\psi\in \Aut L$. Then
 $$ V({\vi})\otimes_{K}V({\psi})\simeq \bigoplus_{c_g\in G/\sim} V(\vi g\psi)^{s_{\vi\psi}^g|c_g|},$$
 where $c_g$ is the equivalence class of $g$, $|c_g|$ its size and
$s_{\vi\psi}^g=\dfrac{|H_{\vi g \psi}|}{|H_{\vi}|
|H_{\psi}|}$.
\end{theorem}

\begin{proof}
Let $\vi, \psi\in \Aut L$. Then by Theorem \ref{proposition-construction-of-balanced}, \eqref{enum-constructed-is-balanced} and Remark \ref{remark-bimod-produces-hom-and-tensor-corresponds-composition}, \eqref{enumerate-tensor}
\begin{align*}
&L\otimes_{K} V(\vi) \otimes_{K}V(\psi)
\simeq
\ds\bigoplus_{g\in
G/H_{\vi}}L_{g \, \vi \, \imath}\otimes_{K} V(\psi)
\simeq \\
\nonumber
&\ds\bigoplus_{g\in
G/H_{\vi}}(L_{g \, \vi}\otimes_{L}L) \otimes_{K}V(\psi)
\simeq
\ds\bigoplus_{g\in
G/H_{\vi}}L_{g \, \vi}\otimes_{L}(L \otimes_{K}V(\psi))
\simeq   \\ \nonumber
&\ds\bigoplus_{g\in
G/H_{\vi}} \ds\bigoplus_{g'\in
G/H_{\psi}} L_{ g \,
\vi}\otimes_{L} L_{g' \, \psi\imath}
\simeq
\ds\bigoplus_{g\in
G/H_{\vi}} \ds\bigoplus_{g'\in
G/H_{\psi}} L_{g \, \vi \, g' \,
\psi\imath}.
\end{align*}
 Then by Corollary \ref{corollary-multoplicities-of-simples}

\begin{equation}\label{formula-10}
[V(\vi)\otimes_{K}V(\psi)]=\sum_{\substack{g\in G/H_{\vi}\\g'\in G/H_{\psi}}} \frac{|H_{g\vi g'\psi}|}{|G|}
 [V(g\vi g'\psi)]= \sum_{c_g\in G/\sim}s_{\vi\psi}^g|c_g|[V(\vi g\psi)].
 \end{equation}

 which completes the proof.

 \end{proof}

\begin{corollary}
\label{corollary-product-simple-balanced} Let $\vi,\psi\in
\Aut L$. Then $\dfrac{1}{|G|}b_{\vi}b_{\psi}\in \bZ[\Aut
L]$ and
 $$ V(b_{\vi})\otimes_{K}V(b_{\psi})\simeq V(\dfrac{1}{|G|}b_{\vi}\cdot b_{\psi}).$$
\end{corollary}

\begin{proof} Clearly,
$$\dfrac{1}{|G|}b_{\vi}b_{\psi}=\sum_{g_{1},g_{2},g\in G}g_{1}\vi g\psi g_{2},$$
which proves the first statement. On the other hand we have
the following equalities in $\bQ[\Aut L]$:
\begin{equation*}
 b_{\vi}\cdot b_{\psi}=
 \big(\sum_{\substack{g\in G/H_{\vi},\\g'\in G}}g\vi g'\big)\big(\sum_{\substack{g\in G/H_{\psi},\\g'\in G}}g\psi g'\big)=
\frac{|G|}{|H_{\vi}||H_{\psi}|}\sum_{g\in G}|H_{\vi
g\psi}|\,b_{\vi g \psi}.
 \end{equation*}
Comparison with \eqref{formula-10} we complete the proof.

\end{proof}

\begin{corollary}
\label{corollary-rational-Grothendieck-isomorphic-Hecke}
The map
$$\Psi:\bQ\otimes_{\bZ}K_{0}(K,L)\rightarrow
\cH_{\bQ}(\Aut L;G),\ \Psi([V(\vi)])=\dfrac{1}{|G|}b_{\vi}$$
 is an isomorphism
of $\bQ$-algebras.
\end{corollary}

\begin{proof}
Since the classes $[V(\vi)]$ and the elements
$\dfrac{1}{|G|}b_{\vi}$, $\vi\in G\backslash \Aut L/G$,
form the $\bQ$-bases in $\bQ\otimes_{\bZ}K_{0}(K,L)$ and in
$\cH_{\bQ}(\Aut L;G)$ respectively, then $\Psi$ is an
isomorphism of $\bQ$-vector spaces. The fact that $\Psi$ is
an algebra homomorphism follows immediately from
Corollary~\ref{corollary-product-simple-balanced}.
\end{proof}

\section{Galois rings}
\label{section-Galois-algebras}

\subsection{Notation and some examples}
\label{subsection-Notation-and-some-examples}

\textit{For the rest of the paper we will assume that $\Ga$ is an
integral
 domain, $K$ the field of fractions of $\Ga$,
$K\subset L$ is a finite Galois extension with the Galois group $G$,
$\imath: K\rightarrow L$ is a natural embedding, $\cM\subset \Aut L$
is a separating monoid on which $G$ acts by conjugations,
$\bar{\Ga}$ is the integral closure of $\Ga$ in $L$, $\cK=(L*\cM)^G$.   }

Recall from
the introduction that an associative  ring $U\subset \cK$ containing $\Ga$ is
called a \emph{Galois $\Ga$-ring}{} if it is finitely generated
over $\Ga$   and
$KU=\cK,\,UK=\cK$. Note that following Lemma
\ref{lemma-Ka-viK-and-right-left-basis} below both equalities in
this definition are equivalent.

\begin{example}
\begin{itemize}
\item  Let $U=\Ga[x;\sigma]$ be the \emph{skew
    polynomial ring} over $\Ga$, where $\sigma\in \Aut
    \Ga$, $x\gamma=\sigma(\gamma)x$, for all $\gamma\in
    \Ga$. Denote $$\cM=\{\sigma^n\mid n=0, 1,
    \ldots\}\subset \Aut K, \cM\simeq \bZ_{+}.$$ Then for $L=K,$ $ G=\{e\}$ the algebra $U$ is a Galois
    $\Ga$-ring in $K*\cM$, when $x$ is identified with $1*\sigma\in K*\cM$.

\item Analogously the \emph{skew Laurent polynomial ring}
    $U=\Ga[x;\sigma^{\pm1}]$ is a Galois ring with
    $\cM=\{\sigma^{n}\,|\, n\in\bZ \}$ and trivial
    $G$.
    
\item Let $\Ga=\k[x_{1},\dots,x_{n}]$ and $\sigma_{1},\dots,\sigma_{n}\in\Aut \Ga$, such that $\sigma_{i}\sigma_{j}=\sigma_{j}\sigma_{i},$ $i,j=$ $1,\dots,n$, $\cM\subset \Aut \Ga$ subgroup generated by $\sigma_{1},\dots,\sigma_{n}.$ Then the skew group ring $\Ga*\cM$  is a Galois ring over $\Ga$ with trivial $G.$

\end{itemize}
\end{example}
More examples will be given in Section \ref{section-Examples-of-Galois-algebras}.

\subsection{Characterization of a Galois ring}
\label{subsection-Definition-of-Galois-algebra} A
$\Gamma$-subbimodule of $\cK$ which for every $m\in
\cM$ contains $[b_{1}m],\dots,[b_{k}m]$ where
$b_{1},\dots,b_{k}$ is a $K$-basis in $L^{H_{m}}$ will be
called a \emph{$\Gamma$-form}{} of $\cK$. We will show
that any Galois subring in $\cK$ is its
$\Gamma$-form.

\begin{lemma}\label{l81}
\label{lemma-Ka-viK-and-right-left-basis}
\label{lemma-existence-of-support-in-a-product}
 Let $u\in U$ be a nonzero element, $T=\supp u$, $u=\ds\sum_{m\in
T/G} [a_m m]$. Then
$$K(\Ga u \Ga)=(\Ga u \Ga)K= K uK\simeq\bigoplus_{m\in T/G} V(m).$$
In particular $U$ is a $\Ga$-form of $\cK$. Besides,
$$L(\Ga u \Ga)=(\Ga u \Ga)L= L uL=\sum_{m\in T} Lm\subset
L*\cM.$$
\end{lemma}

\begin{proof}
Note  that by
Theorem~\ref{proposition-construction-of-balanced}, \eqref{enumerate-K-vi-K} and
Lemma \ref{lemma-separ-monoid}, \eqref{enum-def-of-sep-action-galois}
the modules $V(m),$ $ m\in T/G$, are pairwise non-isomorphic
simple $K$-bimodules.  Since by Lemma
\ref{lemma-separating-action-and-stabilizer},
\eqref{enumerate-H-vi-stabilize-K-vi-K}
$$K[m]K= KK^{m}[m]\simeq V(m), m\in T/G,$$ we have

\begin{align*}
&KuK\subset \sum_{m\in T/G}K [a_m m]K=
\bigoplus_{m\in T/G}K [a_m m]K\simeq \\ & \bigoplus_{m\in T/G}K [ m]K\simeq\bigoplus_{m\in T/G} V(m).
\end{align*}

Since all $V(m)$ are  simple, then the image of $KuK$ in
$W=\osu{m\in T/G}{} V(m)$ generates $W$ as a
$K$-bimodule. Hence $KuK\simeq W$ and therefore $K [a_m
m]K\subset KuK$ for all $m\in T/G$.

For the rest of the proof
it is enough to consider $u=[am]$. Then $$\Ga[am]\Ga=
[\Ga\cdot m(\Ga)a m]\textrm{ and } K\Ga m(\Ga) =K m(K).$$ The
statement $K(\Ga u \Ga)=(\Ga u \Ga)K= K uK$ now follows from Lemma
\ref{lemma-separating-action-and-stabilizer},
\eqref{enumerate-H-vi-stabilize-K-vi-K}.

Obviously $L[am]$
is a $L$-subbimodule in $\su{m\in T}{} Lm$. Since this  is
a direct sum of non-isomorphic simple $L$-bimodules, any
subbimodule has the form $\su{m\in T'}{} Lm, T'\subset T.$
On the other hand $\supp  [am]=T$, and thus
$L[am]=\su{m\in T}{} Lm$.
\end{proof}

\begin{corollary}
\label{corollary-support-of-an-inserted-product}
Let $[a\vi],[b\psi]\in \cK$. Then
$$\supp\, [a\vi]\Ga[b\psi]=\supp\, [a\vi]\supp\,[b\psi]=\cO_{\vi}\cO_{\psi}. $$
\end{corollary}

\begin{proof}
Multiplication on $L$ does not change the support. Then applying Lemma \ref{lemma-existence-of-support-in-a-product}
\begin{align*}
&\supp\, [a\vi]\Ga[b\psi]=\supp L( [a\vi]\Ga[b\psi])=\supp L(K [a\vi]\Ga)[b\psi]=\\
&\supp (L[a\vi]L)[b\psi]=\supp \big(\su{m\in \cO_{\vi}}{} Lm\big)[b\psi]=\cO_{\vi}\cO_{\psi}.
\end{align*}
\end{proof}

\begin{proposition}
\label{proposition-characterization-Galois-by-generators}
Assume  a ring $U\subset \cK$ is generated
over $\Ga$ by  $u_{1},\dots,$ $u_{k}\in U$.
\begin{enumerate}
\item\label{enumerate-support-generates-cM-U-is-Galois}
If\, $\ds\bigcup_{i=1}^{k}\supp u_{i}$
generate $\cM$ as a semigroup, then $U$ is a Galois
ring.

\item\label{enumerate-by-mult-on-L-get-all-is-Galois}
If $LU=L*\cM$, then $U$ is a Galois ring.
\end{enumerate}
\end{proposition}

\begin{proof}
The statement \eqref{enumerate-by-mult-on-L-get-all-is-Galois} follows from \eqref{enumerate-support-generates-cM-U-is-Galois}.
Consider a $K$-subbimodule $Ku_{1}K+$ $\dots+$ $Ku_{k}K$ in
$\cK$. By Lemma
\ref{lemma-Ka-viK-and-right-left-basis}, this bimodule
contains the elements
$[a_{1}\vi_{1}],$ $\dots,$ $[a_{N}\vi_{N}]$, where
$\vi_{1},$ $\dots,$ $\vi_{N}$ is a set of generators of $\cM$.
By Corollary \ref{corollary-support-of-an-inserted-product} $\supp\big([a_1m_1]\Ga[a_2m_2]\big)=\supp[a_1m_1]\cdot
\supp[a_2m_2]$ for $[a_1m_1],[a_2m_2]\in U$, then for
every $m\in \cM$ there exists a nonzero $a_{m}\in
L^{H_{m}}$ such that $[a_{m}m]\in U.$
\end{proof}

\begin{theorem}
\label{theorem-shat-algebras} Let $U$ be a Galois ring,  $e\in \cM$  the unit element  and
$U_{e}= U\cap Le$. Then
\begin{enumerate}

\item \label{enum-GA-is-G-invariant} For every $x\in U$
    holds $x_{e}\in K$ and $U_{e}\subset Ke$.

\item \label{enum-shat-is-dense} The $\k$-subalgebra in
    $L*\cM$ generated by $U$ and $L$ coincides with
    $L*\cM$.

\item \label{enum-shat-is-maximal} $U\cap K$ is a
    maximal commutative $\k$-subalgebra in $U$.

\item \label{enum-shat-center} The center $Z(U)$ of
    algebra $U$ equals $U\cap K^{\cM}$.
\end{enumerate}
\end{theorem}

\begin{proof} Let $x\in U$ and
$x_{e}=\lambda$, $\lambda\in L$. Then for any $g\in G$
holds $\lambda =x_{e}= (x^{g})_{e}=\lambda^{g}$. Hence
$\lambda\in L^{G}=K$. The statement \eqref{enum-shat-is-dense} follows from Lemma \ref{lemma-Ka-viK-and-right-left-basis}.

Consider any $x\in L*\cM$
such that $x\gamma=\gamma x$ for all $\gamma\in \Gamma$.
Assume $x_{\vi}\ne0$ for some $\vi\in\cM,$ $\vi\ne e$.  Since the action of $\cM$ is separating, there  exists
$\gamma\in \Gamma$ such that $\gamma^\vi \ne \gamma$. Then
$(\gamma x)_{\vi}= \gamma x_{\vi}\ne
\gamma^{\vi}x_{\vi}=(x\gamma)_{\vi}$ which is a contradiction.
Hence $x\in U\cap Le=U_e\subset K$ which completes the
proof of \eqref{enum-shat-is-maximal}.

To prove
\eqref{enum-shat-center} consider a nonzero $z\in Z(U)$.
It follows from the proof of \eqref{enum-shat-is-maximal} that $z\in U\cap K$. Besides, $z\in \Ga\cap Z(U)$ if and only if for every
$[a\vi]\in U$ holds $z[a\vi]=[a\vi]z$, i.e. $z=z^{\vi}$.
\end{proof}

Theorem \ref{theorem-shat-algebras}, \eqref{enum-shat-is-maximal} in
particular shows that  an noncommutative associative algebra is never a Galois
ring over its center. For the same reason the universal enveloping
algebra of a simple finite-dimensional Lie algebra is not a Galois
ring over the enveloping algebra of its Cartan subalgebra.

 A submonoid $H$ of $\cM$ is called an \emph{ideal} of $\cM$ if $\cM
H\subset H$ and $H\cM\subset H$.

\begin{corollary}
\label{corollary-ideals-in-K-localization} There is one-to-one
correspondence between  the two-sided ideals in $\cK$ and
the $G$-invariant ideals in the monoid $\cM$. This correspondence
is given by the following bijection \begin{equation}
\label{equation-correspondence} I\longmapsto
\cI=\cI(I)=\bigcup_{u\in I}\supp u,\quad \cI\longmapsto
I=I(\cI)=\sum_{\vi\in\cI} K[\vi] K,
\end{equation}
where $I\subset \cK$, $\cI\subset\cM$ are ideals, $I\ne0$,
$\cI$ is $G$-invariant. In particular, if $\cM$ is a group then
$\cK$ is a simple ring.
\end{corollary}

\begin{proof}
Let $I$ be a nonzero ideal in $\cK$. If $0\ne u\in I$ then
$$KuK\simeq\ds\sum_{\vi\in\supp u/G} K[\vi]K $$ by  Lemma~\ref{lemma-existence-of-support-in-a-product} and for every $m\in\cM$ holds $(K[m]K)(KuK)\subset
I$, $(KuK)(K[m]K)\subset
I$. By Corollary
\ref{corollary-support-of-an-inserted-product} for
every $m\in\cM$ and $\vi\in \supp u$ there exist $u',u''\in I$ such
that $m\vi\in \supp u'$ and  $\vi
m\in \supp u''$. This gives the map $I\mapsto \cI(I)$.
Analogously, $I(\cI)$ is a two-sided ideal in $\cK$ and both
maps are mutually inverse.
\end{proof}

\begin{proposition}
\label{proposition-Ore-condition-for-Gamma} Let $U$ be a Galois
ring  over   $\Ga$, $S=\Ga\setminus \{0\}$.
\begin{enumerate}
\item\label{item-Ore-condition-for-Gamma} The
multiplicative set $S$ satisfies both left and right
Ore condition.

    \item\label{item-U-is-order} The canonical embedding $U\hookrightarrow \cK$ induced the isomorphisms of rings of fractions $[S^{-1}]U\simeq \cK$, $
    U[S^{-1}]\simeq \cK$.
\end{enumerate}
\end{proposition}

\begin{proof}
Assume $s\in S, u\in U$. Following Lemma
\ref{lemma-Ka-viK-and-right-left-basis}, $U$ contains a
right $K$-basis $u_{1},\dots,u_{k}$ of $KuK$, hence in $\cK$ holds
$$s^{-1}u=\ds\sum_{i=1}^{k}u_{i} \gamma_{i} s_{i}^{-1}\textrm{ for
some }s_{i}\in S,\gamma_{i}\in\Ga,i=1,\dots,k.$$  Then in
$U$ holds $$u\cdot (s_{1}\dots
s_{k})=s\cdot(\ds\sum_{i=1}^{k} u_{i}\gamma_{i} s_{1}\dots
s_{i-1}s_{i+1}\dots s_{k}),$$ which shows
\eqref{item-Ore-condition-for-Gamma}.
Besides $S$ acts on $U$ torsion free both from the left and
from the right.
Then there exist the right and left rings
of fractions $U[S^{-1}]$,  $[S^{-1}]U$.
Following Lemma
\ref{lemma-Ka-viK-and-right-left-basis}, the canonical
embedding $U\hookrightarrow \cK$  satisfies the
conditions for the ring of fractions ((i),(ii), (iii), \cite{mcr:nnr}, 2.1.3). Hence
\eqref{item-U-is-order} follows.
\end{proof}

\begin{theorem}
\label{theorem-tensor-product} The tensor product of two Galois
rings is a Galois ring.
\end{theorem}

\begin{proof}
Let $U_i$ be   a Galois $\Ga_i$-subring in the
skew-group algebra $L_i*\cM_i$ with  fraction fields
$K_{i}$, $G_{i}=G(L_{i}/K_{i})$ $i=1,2$. Then
$\cM=\cM_1\times \cM_2$ acts on $L_1\otimes_{\k} L_2$,
$(m_{1},m_{2})\cdot (l_{1}\otimes
l_{2})=m_{1}l_{1}\otimes m_{2}l_{2}$. Since $\k$ is
algebraically closed, $L_1\otimes_{\k} L_2$ is a domain,
hence $\cM$ acts on its  field of fractions $L$. Let $K\subset
L$ be the field of fractions of $K_{1}\otimes_{\k}K_{2}$. The extension $K\subset L$ is
a finite Galois extension with the Galois group
$G=G_{1}\times G_{2}$.
 Consider the composition
$$
\imath:U_{1}\otimes_{\k} U_{2}\myto  (L_{1}*\cM_{1})\otimes_{\k} (L_{2}*\cM_{2})\stackrel{\Phi}\simeq (L_{1}\otimes_{\k}L_{2})*(\cM_{1}\times \cM_{2})\hookrightarrow L*\cM.
$$

We identify $U_{1}\otimes_{\k}U_{2}$ with its image. To endow $U_{1}\otimes_{\k}U_{2}$ with the structure of a Galois
ring we shall prove that $L (U_{1}\otimes_{\k}U_{2})= L*\cM$ (Proposition \ref{proposition-characterization-Galois-by-generators}).
But  $L (U_{1}\otimes_{\k}U_{2})\supset L_{1}U_{1}\otimes_{\k}L_{2}U_{2}=(L_{1}*\cM_{1})\otimes_{\k}(L_{2}*\cM_{2})$, which contains $\Phi^{-1}(\cM_{1}\times \cM_{2})$.
\end{proof}

\section{Galois orders}
\label{section-Integral-Galois-algebra}
\subsection{Characterization of  Galois orders}
\label{subsection-Definition-of-integral-algebras} Let $M$ be a
right $\Gamma$-submodule in a torsion free right $\Ga$-module $N$.
Consider the right
subbimodule in $N$
$$\bD_{r,N}(M)=\{x\in N\,| \text{ there exists } \gamma\in \Gamma,
\gamma\ne 0 \text{ such that } x \cdot\gamma\in M\},$$ which is clearly a right $\Ga$-module. For the left modules $M\subset N$ analogously is defined $\bD_{l,N}(M)$. If $N$ is a
Galois order $U$ over $\Ga$, the we write  $\bD_{r}(M)$ and $\bD_{r}(M)$.

\begin{lemma}
\label{lemma-splitting-of-divisor-closed-submodule} For
right $\Ga$-submodules of $U$ holds the following:
\begin{enumerate}
 \item \label{enum-right-d-idempotentness-of} $M\subset
     \bD_{r}(M)$, $\bD_{r}(\bD_{r}(M))=\bD_{r}(M)$.

\item \label{enum-right-d-and-fraction-field}
    $\bD_{r}(M)=M K \cap U$.

\item\label{enum-right-d-monotonity} If $N\subset M$
    then $\bD_{r}(N)\subset\bD_r(M)$.

\item\label{enum-right-d-for-gamma}
    $\bD_{r}(\Ga)=U_{e}$.

\end{enumerate}
\end{lemma}

\begin{proof} Statements
\eqref{enum-right-d-idempotentness-of} and
\eqref{enum-right-d-monotonity} are obvious. Statement
\eqref{enum-right-d-and-fraction-field}  follows from the fact
that $U$ is torsion free over $\Gamma$. Theorem
\ref{theorem-shat-algebras} \eqref{enum-GA-is-G-invariant}
claims that $U_{e}\subset K$, implying
\eqref{enum-right-d-for-gamma}.
\end{proof}

Lemma \ref{lemma-splitting-of-divisor-closed-submodule}, \eqref{enum-right-d-and-fraction-field} gives the following characterization of  Galois orders (cf. Definition \ref{definition-of-integral}) .

\begin{corollary}
\label{corollary-integral-Galois-algebra} A Galois ring $U$ over
a noetherian $\Ga$ is \emph{ right Galois order}{} if and
only if  for  every finitely generated right $\Gamma$-module     $M\subset U$, the right $\Gamma$-module  $\bD_{r}(M)$ is finitely generated.
\end{corollary}

\begin{corollary}
\label{corollary-integral-is-projective}
If a Galois ring $U$ over a noetherian domain $\Ga$  is projective as a right (left) $\Ga$-module then   $U$ is a right  (left) Galois order.
\end{corollary}

\begin{proof}
 If $U$ is right projective, then there exists some projective right $\Ga$-module $U'$, such that $U\oplus U'\simeq \osu{\cI}{}\Ga$ for some set $\cI$. If $M$ is a finitely generated right submodule in $U$, then there exists a finite subset $\cJ\subset \cI$, such that $M\subset \osu{\cJ}{}\Ga\subset \osu{\cI}{}\Ga$. Then $D_{r,U}(M)= D_{r,U\oplus U'}(M)=D_{r, \osu{\cJ}{}\Ga}(M)\subset \osu{\cJ}{}\Ga$. Then $D_{r}(M)$ is finitely generated since $|\cJ|<\infty$ and $\Ga$ is  noetherian.
\end{proof}

\begin{corollary}\label{Ue-normal}
If $U$ is right (left)  Galois order then $\Ga\subset U_e$
is an integral extension. In particular $U_e$ is a normal ring.
\end{corollary}

\begin{proof}
Lemma~\ref{lemma-splitting-of-divisor-closed-submodule},
\eqref{enum-right-d-for-gamma} shows that $U_e=D_{r}(\Ga)\subset K$
 is finitely generated right (left) $\Ga$-module. Moreover,
it is finitely generated as left and right $\Gamma$-module
simultaneously. The statement  follows from Proposition
\ref{proposition-integrity}.
\end{proof}

 We will show in Theorem
\ref{theorem-for-involutive-algebra-right-is-left},
\eqref{enum-involutive-U-0-is-integral} that the converse statement holds when $\cM$ is a group.

\subsection{Harish-Chandra subalgebras}
\label{subsection-Harish-Chandra-subalgebras}
Following \cite{dfo:hc} a commutative subalgebra ${\Gamma}\subset U$ is
called a \emph{Harish-Chandra subalgebra}{} in $U$  if for
any $u\in U$, the ${\Gamma}-$bimodule $\Ga u\Ga$ is
finitely generated both as a left and as a right
$\Ga$-module. Assume $\Gamma$ and some
family  $\{u_{i}\in U\}_{i\in I}$ generates $U$ as
$\k$-algebra and every $\Gamma u_{i}\Gamma,i\in I$ is left
and right finitely generated. Then it is easy to see, that
$\Gamma$ is a Harish-Chandra subalgebra in $U$.

\begin{proposition}
\label{proposition-integral-Harish-Chandra} Assume  $\Ga$ is
finitely generated algebra over $\k$, $U$ is a Galois
ring. Then $\Ga$ is a Harish-Chandra subalgebra in $U$ if and only
if $m\cdot \bar{\Ga}= \bar{\Ga}$ for every $m\in\cM$.
\end{proposition}

\begin{proof}
Note that $\bar{\Ga}$ is finitely generated as
$\Ga$-module (Proposition \ref{proposition-integrity}). Suppose first $m\cdot \bar{\Ga}= \bar{\Ga}$
for every $m\in \cM$.
  It is enough to
prove that $\Ga[a m ]\Ga$ is finitely generated as a left
(right) $\Ga$-module for any  $m\in \cM,a\in L$. But following \eqref{space-structure}
\begin{equation}
\label{equation-moving-ga-and-vi} \Ga [a m ] \Ga=[\Ga\cdot
 m (\Ga)a m ]=[a m \Ga\cdot m ^{-1}(\Ga)]\end{equation} is
finitely generated over $\Ga$ from the left, since $\Ga m
(\Ga)\subset\bar{\Ga}$, and it is finitely generated from
the right, since $\Ga m ^{-1}(\Ga)\subset \bar{\Ga}$.
Conversely, assume $\Ga [a m ]\Ga$ is finitely generated
right $\Ga$-module for any  $[a m ]\in U$. It means
that $\Ga\cdot
 m ^{-1}(\Ga)$ is finite over $\Ga$, i.e.
$ m ^{-1}(\Ga)\subset\bar{\Ga}$. Analogously, $ m
(\Ga)\subset \bar{\Ga}$.
\end{proof}

\begin{proposition}\label{integral-hchandra}
If $\, U$ is a right (left)  Galois order over a noetherian
$\Ga$ then for any $m\in \cM$ holds $m^{-1}(\Ga)\subset \bar{\Ga}$
($m(\Ga)\subset \bar{\Ga}$).
\end{proposition}

\begin{proof}
Let $U$ be right  Galois order, $[a m ]\in U$,  $\ga\in\Ga$. Assume $x= m
^{-1}(\gamma)\not\in\bar{\Ga}$. Then
 the right $\Ga$-submodule of $U$,
$$M=\ds\sum_{i=0}^{\infty}\ga^{i}[a m ]\Ga=
\sum_{i=0}^{\infty}[a m x^{i}\Ga] $$ is not finitely generated.
On the other hand, $x$ is an algebraic element over $K$. Let
$$\ga_{0}x^{n}+\ga_{1}x^{n-1}+\dots+\ga_{n}=0,\
\ga_{i}\in\Ga,\ \ga_{0}\ne0.$$ Consider the following
finitely generated
 right $\Ga$-module
$N=\su{i=0}{n-1}\ga^{i}[a m ]\Ga=
\su{i=0}{n-1}[a m  x^{i}\Ga].$ But
$M\subset\bD_{r}(N)$ which is a contradiction.  The case of
left order treated analogously.
\end{proof}

From Proposition~\ref{integral-hchandra} and
Proposition~\ref{proposition-integral-Harish-Chandra} we
immediately obtain

\begin{corollary}
\label{corollary-right-left-integra-is-hc} Let $\Ga$ be a noetherian
domain and $U$  a Galois order over   $\Ga$. Then
$\Ga$ is a Harish-Chandra subalgebra in $U$.
\end{corollary}

\begin{remark}\label{rem-61} Let $\Ga$ be integrally closed in $K$ and $\vi:K\myto K$   an automorphism of infinite order, such that
$\vi(\Ga)\nad{\subset}{\ne}\Ga$. Set $L=K$,
$\cM=\{\vi^{n}|n\geq 0\}$. Then $L*\cM$ is isomorphic to
the skew polynomial algebra $K[x;\vi]$ (\cite{mcr:nnr}).
Its subalgebra $U$ generated by $\Ga$ and $x$ is a Galois
ring. Clearly, $U$ is  left Galois order (but not right Galois order).
\end{remark}

\subsection{Properties of  Galois orders}
\label{subsection-Properties-of-Galois-algebras} Let $U$ be a Galois
ring over   $\Ga$, $S\subset \cM$  a finite $G$-invariant subset.
Denote \begin{equation}
\label{equation-definition-of-Us}U(S)=\{u\in U\, |\,\supp u\subset S\}.
\end{equation}
Obviously, it is a $\Ga$-subbimodule in $U$ and $\bD_{r}(U(S))=\bD_{l}(U(S))=U(S)$.
This notion will give us one more characterization of Galois orders (Theorem \ref{theorem-integrality-equivalence-of-defs}).

 It will be convenient to
consider the $\Ga$-bimodule structure of $U$ as a
$\Ga\otimes_{\k}\Ga$-module structure. For every $f\in
\Gamma$ define $f_{S}^r\in \Gamma\otimes_{\k}L$
(respectively $f_{S}^l\in L\otimes_{\k}\Ga$) as follows
\begin{align}
\label{equation-relations-in-as-in-gz}
&f_{S}^r=\prod_{s\in S}(f\otimes 1 -
1\otimes
f^{s^{-1}})=\sum_{i=0}^{|S|}
f^{|S|-i}\otimes h_{i}, \
h_{0}=1, \\
& (\textrm{resp.} f_{S}^l=\prod_{s\in S}(f^s\otimes
1 - 1\otimes f)=\sum_{i=0}^{|S|}
h_{i}'\otimes f^{|S|-i}, \
h'_{0}=1).
\end{align}

Since $S$ is $G$-invariant, then all $h_{i}$ and $h_{i}'$ are $G$-invariant expressions in $f^{m},$ $m\in \cM$, they belongs to $K$.
If $U$ is right (left) integral, then $h_{S}^{r}\in \Ga\otimes\Ga$ ($h_{S}^{l}\in \Ga\otimes\Ga$).
We will consider the properties of
$f_S=f_{S}^r$, the case of $f_{S}^l$ can be treated
analogously.

\begin{lemma}
\label{lemma-relation-as-in-gz-acts-zero} Let
$m^{-1}(\Ga)\subset \bar{\Ga}$ for any $m\in \cM$,
$S\subset \cM$  a $G$-invariant subset, $u\in U$,  $f\in
\Gamma$.
\begin{enumerate}

\item\label{enum-f-S-kills-U-S}  $u\in U(S)$ if and
    only if $f_{S}\cdot u=0$ for every $f\in\Gamma$.

\item \label{enum-f-complement-to-S-throws-U-S} If
    $T=\supp u \setminus S$ then $f_{T}\cdot u\in
    U(S)$ for every $f\in\Gamma$.

\item\label{enum-f-S-mult-on-standard} If
    $f_{S}=\su{i=1}{n}f_{i}\otimes g_{i}$,
    $[am]\in L*\cM$ then $f_{S}\cdot
    [am]=[(\su{i=1}{n} f_{i}g_{i}^{m}a)m]$ $=[\pru{s\in S }{}(f -
    f^{ms^{-1}})am]$.

\item\label{enum-f-S-mult-components-homomorphism} Let
    $S$ be  a $G$-orbit and $T$ an $G$-invariant subset
    in $\cM$. The $\Ga$-bimodule homomorphism
    $P^{T}_{S}(=P^{T}_{S}(f)):U(T)\myto U(S)$,
    $u\mapsto f_{T\setminus S}\cdot u,f\in \Gamma$ is
    either zero or $\Ker P_{S}^{T}=U(T\setminus S)$
    (both cases are possible, cf.
    \eqref{enum-f-S-kills-U-S}).

\item\label{enum-f-S-mult-split-components-mono} Let
    $S=S_{1}\sqcup\dots\sqcup S_{n}$ be the
    decomposition of $S$ in $G$-orbits and
    $P_{S_{i}}^{S}:U(S)\myto U(S_{i})$ for some
    $f_{i}\in \Gamma$, $i=1,\dots,n$ are defined in
    \eqref{enum-f-S-mult-components-homomorphism}
    nonzero homomorphisms.  Then the homomorphism
\begin{equation}
\label{equation-bimodule-distinguishing-components}
P^{S}:U(S)\myto\bigoplus_{i=1}^{n}
U(S_{i}),\
P^{S}=(P_{S_{1}}^{S},\dots,P_{S_{n}}^{S}),
\end{equation}
is a monomorphism.

\item\label{enum-f-S-holds-for-hc}
The statements above hold, if $\Ga$ is a Harish-Chandra subalgebra in $U$.
\end{enumerate}
\end{lemma}

\begin{proof}
Consider any $[am]\in L*\cM$, $s\in \Aut L$. Then
$$(f\otimes 1 - 1\otimes f^{s})\cdot
[am]=[fam]-[amf^{s}]= [(f-f^{ms})am]\text{, hence}$$
\begin{align*}
&f_{S}\cdot [am]=\prod_{s\in S}(f\otimes 1 - 1\otimes
f^{s^{-1}})\cdot[am] = [\prod_{s\in S }(f -
f^{ms^{-1}})am].\end{align*} If $m\in S$, then one of
$f - f^{ms^{-1}}$ equals zero, hence, $f_{S}\cdot [am]=0$. To
prove  the converse we show that for any $m\not\in S$ there
exists $f\in\Ga$ such that $f\ne f^{ms^{-1}}$ for all $s\in
S$. Following Lemma \ref{lemma-separ-monoid}, \eqref{enum-def-of-sep-action-act} for
every $m\in\cM, m\ne e$, the space of $m$-invariants
$\Ga^{m}\ne\Ga$. But the $\k$-vector space $\Ga$ can not be
covered by  finitely many proper subspaces $\Ga^{ms^{-1}},
{s\in S}$, that completes the proof of
\eqref{enum-f-S-kills-U-S}. Obviously, $f_{\supp u}\cdot u=0$ for any $f\in \Gamma$.
Then  statement \eqref{enum-f-complement-to-S-throws-U-S}
follows from \eqref{enum-f-S-kills-U-S} and from the
equality $f_{\supp u}=f_{S}f_{T}$. Statement
\eqref{enum-f-S-mult-on-standard} follows from the formulas
\eqref{align-def-a-vi-b}, \ref{subsection-Skew-group-rings-and-skew-group-categories}.
By \eqref{enum-f-S-mult-on-standard}, $f_{T\setminus
S}\ne0$ if and only if $\underset{i=1}{\overset{n}\sum}
f_{i}g_{i}^{m}\ne0$, and in this case $f_{T\setminus
S}$ acts on $U(S)$ injectively, that proves \eqref{enum-f-S-mult-components-homomorphism}.

Finally, 
\eqref{enum-f-S-mult-split-components-mono} follows from
\eqref{enum-f-S-mult-components-homomorphism}, since
$\underset{i=1}{\overset{n}\cap}\Ker P_{S_{i}}^{S}=0$ and \eqref{enum-f-S-holds-for-hc} follows from the definition.
\end{proof}

\begin{theorem}
\label{theorem-integrality-equivalence-of-defs} Let $U$ be a Galois
ring over  a noetherian
 Harish-Chandra subalgebra $\Ga$. Then the following
statements are equivalent:

\begin{enumerate}

\item\label{enum-equiv-of-defs-int-right} $U$ is right
    (respectively left) Galois order.

\item\label{enum-equiv-of-defs-int-right-any-s} $U(S)$
    is finitely generated right (respectively left)
    $\Ga$-module for any finite $G$-invariant $S\subset
    \cM$.

\item\label{enum-equiv-of-defs-int-right-elements}
    $U(\cO_{m})$ is finitely generated right
    (respectively left) $\Ga$-module for any $m\in\cM$.

\end{enumerate}

\end{theorem}

\begin{proof} Assume $U$ is right
 Galois order. Let $S$ be a finite $G$-invariant  subset of $\cM$,  and  $u_{1},$ $\dots,$ $u_{k}\in $ $U(S)$ a basis of $ U(S)K$ as a  right $K$-space. Then
$$ \bD_{r}\big(\ds\sum_{i=1} ^{k}
u_{i}\Gamma \big)=\big(\ds\sum_{i=1} ^{k} u_{i}\Gamma \big)K\cap
U= U(S)K\cap U= \bD_{r}(U(S))=U(S).$$ Therefore,
$U(S)=\bD_{r}\big(\ds\sum_{i=1} ^{k} u_{i}\Gamma\big)$,
which proves \eqref{enum-equiv-of-defs-int-right-any-s}.
Obviously, \eqref{enum-equiv-of-defs-int-right-any-s}
implies \eqref{enum-equiv-of-defs-int-right-elements}.
Assume \eqref{enum-equiv-of-defs-int-right-elements} holds. Let $M\subset U$ be a
finitely generated right $\Ga$-submodule, $S=\supp M$. Then $M\subset U(S)$  and $\bD_r(M)\subset \bD_r(U(S))=U(S)$. By Corollary~\ref{corollary-integral-Galois-algebra},
 it remains to prove that $U(S)$ is
finitely generated. Let
$S=S_{1}\sqcup\dots\sqcup S_{n}$ be the decomposition of
$S$ into $G$-orbits. Following Lemma
\ref{lemma-relation-as-in-gz-acts-zero},
\eqref{enum-f-S-mult-split-components-mono},  $P^S$  embeds $U(S)$ into $\osu{i=1}{n}U(S_{i})$ that completes the proof.
\end{proof}

\begin{theorem}
\label{theorem-for-involutive-algebra-right-is-left} Assume that
$U$ is a Galois ring,   $\Ga$ is noetherian and
$\cM$ is a group.
\begin{enumerate}

\item \label{enum-involutive-right-integral}
Assume $m^{-1}(\Gamma)\subset \bar{\Gamma}$ (resp.
$m(\Gamma)\subset \bar{\Gamma}$). Then $U$ is right (resp. left) Galois order if and only if $\,U_{e}$ is an integral
extension of $\Gamma$.

\item \label{enum-involutive-U-0-is-integral}
Assume $\Gamma$ is a Harish-Chandra subalgebra in $U$. Then $U$ is
 a Galois order if and only if $\,U_{e}$ is an integral
extension of $\Gamma$.
\end{enumerate}
\end{theorem}

\begin{proof} Obviously \eqref{enum-involutive-U-0-is-integral} follows from
\eqref{enum-involutive-right-integral} and Proposition~\ref{proposition-integral-Harish-Chandra}. The statement ``only if'' in \eqref{enum-involutive-right-integral} follows from Corollary \ref{Ue-normal}.  Assume $U_{e}$ is an integral extension of $\Gamma$, $m^{-1}(\Gamma)\subset \bar{\Gamma}$, but $U$ is not right
order. Following Theorem
\ref{theorem-integrality-equivalence-of-defs},
\eqref{enum-equiv-of-defs-int-right-elements} there exists
$m\in\cM$, such that $U(\cO_{ m})$ is not finitely
generated.

Since $\cM$ is a group  then  there exists
$[bm^{-1}]\in U$ by Lemma
\ref{lemma-Ka-viK-and-right-left-basis}.  Since $H_{m}=H_{m^{-1}}$ for any nonzero  $\ga\in\Ga$ holds
\begin{equation}
\label{equation-formula-for-unit-component}
([bm^{-1}] \ga[ma])_{e}=\sum_{g\in G/H_{m}} b^{g} \ga^{(m^{-1})^{g}}a^{g}.
\end{equation}
Denote this expression by $v_{\ga}(a)$, $\ga\in \Ga, a\in L^{H_{m}}$. Then $v_{\ga}: L^{H_{m}}\to K$ is a right $K$-linear map and $v_{\ga_{1}}+v_{\ga_{2}}=v_{\ga_{1}+\ga_{2}},\ga_{1},\ga_{2}\in \Ga$.

Denote $|G/H_{m}|$ by $n$.  Let $\{a_{i}\in L^{H_{m}}\mid i=1,\dots,n\}$ be a basis of $L^{H_{m}}$ over $K$. In particular, $[ma_{i}], i=1,\dots,n$ form a right $K$-basis of $KmK$. It will be convenient for us to enumerate entries of matrices both by the classes from $G/H_{m}$ and the numbers $1,$ $\dots,$ $n.$

\begin{lemma}
\label{lemma-cryterion-integrity-matrices}
\begin{enumerate}
\item
\label{item-cryterion-integrity-ba}
For any $b\in L^{H_{m}},$ $b\ne0$ the $G/H_{m}\times n$ matrix over $L$
 $$X=\big(b^{g} a_{i}^{g}\mid  g\in G/H_{m}; i=1,\dots,n\big)$$ is non-degenerated.

\item

\label{item-cryterion-integrity-ga} There exists $\ga_{1},\dots,\ga_{n}\in \Ga$, such that $n\times G/H_{m}$ matrix
$$ Y=\big(\ga_{i}^{g m^{-1}g^{-1}} \mid   i=1,\dots,n; g\in G/H_{m}\big)$$
is non-degenerated. Besides for $n\times n$ matrices holds
$$YX=\big(v_{\ga_{i}}(a_{j})\mid i,j=1,\dots,n\big).$$

\item
\label{item-cryterion-integrity-dual-basis}
Let $Z=\big(\mu_{ij}\mid i=1,\dots,n; j=1,\dots,n\big)$ be a non-degenerated matrix over $K$, $b_{i}=\su{j=1}{n} a_{j}\mu_{ij},\, i=1,\dots,n $ the new basis of $L^{H_{m}}$. Then
$$ (YX)Z=\big(v_{\ga_{j}}(b_{i})\mid i,j=1,\dots,n\big).$$

\item
\label{item-cryterion-integrity-inverse}
In particular, if $Z=(YX)^{-1}$ holds $$v_{\ga_{i}}(b_{j})=\delta_{ij},i,j=1,\dots,n.$$

\end{enumerate}

\end{lemma}

\begin{proof}
To prove the first statement there is enough to prove the invertibility of the matrix $\big(a_{i}^{g}\mid  g\in G/H_{m}; i=1,\dots,n\big)$. Assume, opposite, i.e. $(\su{g\in G/H_{m}}{}\la_{g}g) (a_{i})=0, \la_{g}\in L$ for some  vector $(\la_{g}\mid g\in G/H_{m})\ne0$ and for any $i=1,\dots,n$.
Then $\big(\su{g\in G/H_{m}}{}\la_{g}g\big)|_{L^{H_{m}}}=0$, which contradicts to the independence of different characters $g|_{L^{H_{m}}}:L^{H_{m}}\to L$, $g\in G/H_{m}.$

Analogously all  $\{{g m^{-1}g^{-1}} \mid  g\in G/H_{m}\}$
act differently in restriction on $\Ga$, hence the row rank of $G/H_{m}\times \Ga$ matrix over $L$ $$\big(\ga^{g}\mid g\in G/H_{m}; \ga\in \Ga
 \big),$$ equals $n$. Then its column rank of this matrix equals $n$ as well, that finishes the proof of the second statement.

The third and fours statement is proved by direct calculation

$$ (YX)_{ij}=\sum_{g\in G/H_{m}} b^{g} \ga_{i}^{g m^{-1}g^{-1}}a_{j}^{g} =v_{\ga_{i}}(a_{j}).$$

$$ ((YX)Z)_{ij}=\sum_{l=1}^{n} v_{\ga_{i}}(a_{l})\mu_{lj}=v_{\ga_{i}}(\sum_{l=1}^{n}a_{l}\mu_{lj}) =v_{\ga_{i}}(b_{j}).$$

The last statement is obvious.
\end{proof}

Assume  $U(\cO_{ m})$ contains a strictly ascending chain of
right $\Gamma$-submodules
\begin{equation}
\label{equation-ascending-chain-right-submodules}
M_{k}=\sum_{i=1}^{k}[m t_{i}]\Ga, \,i=1,2,\dots, \, M=\bigcup_{k=1}^{\infty} M_{k}.
\end{equation}

Fix  $\ga_{1},\dots,\ga_{n}$ from Lemma \ref{lemma-cryterion-integrity-matrices}, \eqref{item-cryterion-integrity-ga} and the  basis $b_{1},\dots,b_{n}$ from Lemma~\ref{lemma-cryterion-integrity-matrices}, \eqref{item-cryterion-integrity-inverse}.


Consider the decomposition $t_{i}=\su{j=1}{n} \ga_{ij}b_{j},\ga_{ij}\in K$. Then there exists $1\leq l\leq n$, such that the $\Ga$-module $T_{l}=\su{i=1}{\infty}\Ga \ga_{il}\subset K$ is not finitely generated. In opposite  case from notherianity of $\Ga$ and $M\subset \bosu{i=1}{n} T_{i}$ follows, that $M$ is finitely generated.

Then   $([b m]\ga_{l}[m^{-1}M])_{e}=v_{\ga_{l}}(M)=T_{l}$, which is not finitely generated.
Let $S=\cO_{m^{-1}}\cO_{m}$. Since $m^{\pm1}(\Gamma)\subset
\bar{\Gamma}$ there exists $F=\ds\sum_{i=1}^{n}
f_{i}\otimes g_{i}\in\Ga\otimes_{\k}\Ga$ (by Lemma
\ref{lemma-relation-as-in-gz-acts-zero},
\eqref{enum-f-S-mult-on-standard}), which defines a nonzero
morphism $P^{S}_{e}:U(S)\myto U(\{e\})=U_e$. Then by Lemma \ref{lemma-relation-as-in-gz-acts-zero}, \eqref{enum-f-S-mult-on-standard}
$$P^{S}_{e}([b m]\ga_{l}[m^{-1}M])= \ga T_{l}\subset U_{e},\, \ga=\ds\sum_{i=1}^{n}
f_{i} g_{i},$$ which means that $U_{e}$ is not finitely generated.
\end{proof}

\begin{corollary}
\label{corol-lazy-is-galois} Let $\cM$ be a
group, $\Gamma$ normal and noetherian, $\cM\cdot
\bar{\Ga}=\bar{\Ga}$,  $\vi_{1},\dots,$ $\vi_{n}\in \cM$  a
set of generators of $\cM$ as a semigroup,
$a_{1},\dots,a_{n}\in \bar{\Ga}$. Then
the subalgebra $U$ in $\cK$ generated by $\Ga$ and
$[a_{1}\vi_{1}],\dots,[a_{n}\vi_{n}]$ is a  Galois
order over  $\Ga$.
\end{corollary}

\begin{proof}
Since $\cM\cdot
\bar{\Ga}=\bar{\Ga}$  any $u\in U$ has a form $u=\underset{m\in
\cM}\sum [a_{m}m],$ where all $a_{m}$ are in $\bar{\Ga}$. In
particular, if $u\in U_{e}$ then $u=[a_{e}e]$ where
$a_{e}\in K\cap \bar{\Ga}$. Since $\Ga$ is normal  $U_{e}=\Ga$. Applying Theorem
\ref{theorem-for-involutive-algebra-right-is-left},
\eqref{enum-involutive-U-0-is-integral} we obtain the
statement.
\end{proof}

The next corollary is a noncommutative analog of Proposition~\ref{prop-nonsingular-commutative-lifting}.

\begin{corollary}
\label{corollary-criterion-integra-algebra} Let $U\subset L*\cM$ be
a Galois ring over noetherian $\Ga$, $\cM$ a group and $\Ga$ a
normal $\k$-algebra. Then the following statements are equivalent
\begin{enumerate}
\item \label{enum-algebra-is-integral} $U$ is a Galois order.

\item \label{enum-divisor-criterion-integrality} $\Ga$ is a Harish-Chandra subalgebra and, if for $u\in U$ there exists a nonzero $\ga\in \Ga$ such that $\ga
    u\in \Ga$ or $u\ga\in \Ga$, then $u\in\Gamma$.
\end{enumerate}
\end{corollary}

\begin{proof} Assume
\eqref{enum-algebra-is-integral}. Then $\Ga$ is a
Harish-Chandra subalgebra  by
Corollary~\ref{corollary-right-left-integra-is-hc}. If $u\gamma\in \Ga$ for  $u\in U$ and $\gamma\in
\Ga$, then $\supp u=\{e\}$, hence $u\in U_{e}$. Applying Corollary~\ref{Ue-normal} we obtain
\eqref{enum-divisor-criterion-integrality}.
To prove the converse
implication consider $u\in U_{e}$. Since $U_{e}\subset K$
(Theorem \ref{theorem-shat-algebras},
\eqref{enum-GA-is-G-invariant}), there exists $\gamma\in
\Gamma$, such that $\gamma u\in \Gamma$. Thus,
$u\in\Gamma$. Theorem
\ref{theorem-for-involutive-algebra-right-is-left},
\eqref{enum-involutive-U-0-is-integral} completes the
proof.
\end{proof}

\subsection{Filtered Galois orders}
\label{subsection-Filtered-Galois-orders}
 Let $U$ be a Galois ring over a noetherian normal
$\k$-algebra $\Gamma$. Suppose in addition that
  $U$ is an
algebra over ${\k}$, endowed with an increasing exhausting
filtration $\{U_{i} \}_{i\in\mathbb Z}$, $U_{-1} =$ $\{ 0
\}$, $U_{0} =$ ${\k}$, $U_iU_j\subset U_{i+j}$ and
$\gr U=$ $ \bosu{i=0}{\infty}
U_{i }/U_{i-1 }$ the associated graded algebra.

 The filtration on $\Gamma$ induces a
 degree "$\deg$" both on $U$ and $\gr{U}$. For $u\in U$
denote by $\bar{u}\in \gr{U}$ the corresponding
homogeneous element and denote by $\gr{\Gamma}$ the image
of $\Ga$ in $\gr{U}$.

\begin{proposition}
\label{prop-reduction-to-commutative} Assume $\gr{U}$ is a domain.  If the canonical
embedding $\imath: \gr{\Gamma} \hookrightarrow \gr{U}$
induces an epimorphism $$\imath^{*}:\Specm \gr{U} \to
\Specm \gr{\Gamma}$$ then $U$ is a Galois order over $\Ga$.
\end{proposition}

\begin{proof}
We apply Corollary
\ref{corollary-criterion-integra-algebra}. Suppose $y=x
u\ne0$, $y,x\in \Ga$, $u\in U\setminus \Ga$ with minimal
possible $\deg y$. Then $\bar{y}=\bar{x} \bar{u}\ne0$ in
$\gr{U}$. By Proposition
\ref{prop-nonsingular-commutative-lifting} $\bar{u}\in
\gr{\Gamma}$.  Hence $\bar{u}=\bar{z}$ for some in $z\in
\Ga$. Since $z\neq u$, we have $y_{1}=x u_{1}$ where
$u_{1}=u-z$, $y_{1}=y-x z$. Then $x,y_{1}\in \Gamma$,
$u\not\in \Gamma$ and $\deg y_{1}<\deg y$. Obtained
contradiction shows that $u\in \Ga$.
\end{proof}

%
%
%
%
%
%
%
%
%
%

\section{Gelfand-Kirillov dimension of Galois orders}
\label{section-Gelfand-Kirillov-dimension-of-Galois-algebras}

In this section we assume that $\cM$ is a group of finite
growth and $\Ga$ is an affine $\k$-algebra of  finite
Gelfand-Kirillov dimension.

\subsection{Growth of group algebras}
Let $S_{*}=\{S_{1}\subset S_{2}\subset\dots \subset
S_{N}\subset \dots\}$ be an increasing chain of finite
sets. Then the growth of $S_{*}$ is defined as
\begin{equation}
\label{equation-definition-of-growth}
\gro(S_{*})=\overline{\lim_{N\to\infty}} \log_{N}|S_{N}|.
\end{equation}
For $s\in S=\ds\bigcup_{i=0}^{\infty}S_{i}$ set $\deg s=i$
if $s\in S_{i}\setminus S_{i-1}$. Let
$\{\ga_{1},\dots,\ga_{k}\}$ be a  set of generators of
$\Ga$. For $N\in\bN$ denote by $\Ga_{N}\subset \Ga$ the
subspace of $\Ga$ generated  by the products
$\ga_{i_{1}}\dots\ga_{i_{t}}$, for all $t\leq N$,
$i_{1},\dots,i_{t}\in \{1,\dots,k\}$. Let
$d_{\Ga}(N)=\dim_{\k} \Ga_{N}$ and let $B_{N}(\Ga)$ be a
basis in $\Ga_{N}$
($B_{1}({\Ga})=\{\ga_{1},\dots,\ga_{k}\}$). Fix a set of
generators of $\cM$ of the form $\cM_{1}=\cO_{\vi_{1}}\cup$
$\dots$ $\cup \cO_{\vi_{n}}$. For $N \geq1$, let $\cM_{N}$
be the set of words $w\in\cM$ such that $l(w)\leq N$, where
$l$ is the length of $w$, i.e.
\begin{align}
\label{align-next-set} \cM_{N+1}=\cM_{N}\bigcup
\bigg(\bigcup_{\vi\in\cM_{1}}\vi\cdot\cM_{N}\bigg).
\end{align}
 Note that all sets
$\cM_{N}$ are $G$-invariant.  Denote the cardinality of
$\cM_{N}$ by $d_{\cM}(N)$.  Let $\cM_{*}=\{\cM_{1}\subset
\cM_{2}\subset\dots \subset \cM_{N}\subset \dots\}$.  Then
$\gro(\cM)$ is by definition $\gro(\cM_{*})$.

Let $\Ga[\cM]$ be the group algebra of $\cM$. Assume, $G$
acts on $\Ga[\cM]$, acting by $\cM$ by conjugations and
trivially on $\Ga$. Then the space $\Ga[\cM]_{N}$ has a
$G$-invariant basis
\begin{equation}
\label{equation-basis-in-the-group-algebra}
B_{N}(\Ga[\cM])=\bigsqcup_{i=0}^{N} \bigsqcup_{\substack{w\in
\cM_{N-i},\\ l(w)=N-i, }}B_{i}(\Ga)w.
\end{equation}
and $\gkdim \Ga[\cM]=\gro{B_{*}(\Ga[\cM])}$. In particular
(e.g. \cite{mcr:nnr}, Lemma 8.2.4)
\begin{equation}
\label{equation-formula-for-GK-dimension-of-group-algebra} \gkdim
\Ga[\cM]=\gkdim \Gamma+\gro(\cM).
\end{equation}

The  growth of the chain $B_{*}(\Ga[\cM])/G$ is equal to
$\gro{B_{*}(\Ga[\cM])}$, since
$$|B_{N}(\Ga[\cM])|>|B_{N}(\Ga[\cM])/G|\geq
\dfrac{1}{|G|}|B_{N}(\Ga[\cM])|.$$

Without loss of generality we can assume that the Galois
ring $U$ is generated over $\Ga$ by a set of generators
$\cG=\{[a_{1}\vi_{1}],$ $\dots,$ $[a_{n}\vi_{n}]\}$. Set
$B_{1}(U)=B_{1}(\Ga)\sqcup \cG$. As above, define the
subspaces $U_{N}$ and dimensions $d_{U}(N)$. For every
$N\geq 1$ fix a basis $B_{N}(U)$  of $U_{N}$.

\subsection{Gelfand-Kirillov dimension}
\label{subsection-Gelfand-Kirillov-dimension}

 The goal of this
section is to prove (under a certain condition) an analogue of the
formula \eqref{equation-formula-for-GK-dimension-of-group-algebra}
for Galois orders.

\begin{condition}
\label{condition-on-the-global-action}  For every finite
dimensional $\k$-vector space $V\subset \bar{\Ga}$ the set
$\cM\cdot V$ is contained in a finite dimensional subspace
of $\bar{\Ga}$.
\end{condition}

\begin{theorem}
\label{theorem-GK-dimension-of-some-Galois} If $U$ is a Galois $\Ga$-ring which satisfies Condition
\ref{condition-on-the-global-action} and $\cM$ is a group  of
finite $\gro(\cM)$, then
\begin{equation}
\label{equation-formula-for-GK-dimension-of-some-Galois} \gkdim
U\geq\gkdim \Gamma+\gro(\cM).
\end{equation}
\end{theorem}

The proof of this result is based on the following lemmas.

\begin{lemma}
\label{lemma-on-an-estimate-of-growth-of-Galois-algebra} If
for some $p,q\in\bZ$ and $C> 0$  for any $ N\in \bN$ holds
\begin{equation}
\label{equation-the-main-unequality-for-GK} d_{U}(pN+q)\geq C
d_{\Ga[\cM]}(N),
\end{equation}
then $\gkdim U\geq \gkdim \Ga[\cM]$.
\end{lemma}

\begin{proof}
\begin{align*}
&\gkdim
\Ga[\cM]=\overline{\lim_{N\to\infty}} \log_{N}d_{\Ga[\cM]}(N)\leq
\overline{\lim_{N\to\infty}} \log_{N}d_{U}(pN+q)=\\
&
\overline{\lim_{N\to\infty}} \log_{pN+q}d_{U}(pN+q)\leq
\overline{\lim_{N\to\infty}} \log_{N}d_{U}(N)=\gkdim U.
\end{align*}
\end{proof}

 Denote by
$N(i)$, $i=1,2,\dots$ the minimal number  such that for
any $m\in\cM_{i}$, $U_{N(i)}$ contains an element of the form
$[bm]$, $b\ne0$.

\begin{lemma}
\label{lemma-main-shift-in-growth}

\begin{enumerate}
\item\label{enumerate-is-in-support} For every
    $i=1,\dots,n$ there exists a finite dimensional
    over $\k$ space $V_{i}\subset \Ga$, such that for
    any $x\in U$ and $m\in \supp x$ there exists $y\in
    [a_{i}\vi_{i}]V_{i}x$ such that $\vi_{i}m\in\supp
    y$. Besides  $|\supp y|\leq |G| |\supp x|$ and
    $\deg y- \deg x\leq d$ for some fixed $d>0$.

\item\label{enumerate-extracted-from-support} For every
    $k\geq1$ there exists $t(k)\geq0$ with the
    following property: for every $j\geq1$ and $u\in
    U_{j}$, such that $|\supp u|\leq k$ and for  any
    $m\in\supp u$ there exists a nonzero element
    $[bm]\in U_{j+t(k)}$.

\item\label{enumerate-globally-presented-in-support}
    The sequence $N(i+1)-N(i)$, $i=1,2,\dots$ is
    bounded.
\end{enumerate}
\end{lemma}

\begin{proof}
Let $L(G/H_{\vi_{i}})$ be the vector space over $L$ with the basis, enumerated by cosets $G/H_{\vi}$, $\vi\in \cM$. We endow this space with the
standard scalar product.  Fix $i$, $1\leq i\leq n$ and consider the  nonzero
vector $$v(x)=( a_{i}^{g} x_{(\vi_{i}^{g})^{-1}\vi_{i}
m}^{\vi_{i}^{g}})_{g\in G/H_{\vi_{i}}}\in L(G/H_{\vi_{i}}).$$ Then
for any $\ga\in \Ga$  immediate calculation shows, that
\begin{align}
\label{align-coefficient}
([a \vi_{i}]\ga x)_{\vi_{i}
m}= v(x)\cdot (\ga^{\vi_{i}^{g}})_{g\in G/H_{\vi_{i}}}\in L^{H_{\vi_{i}m}}.
\end{align}

Since  $\vi_{i}^{g},g\in G/H_{\vi}$ are different,  there
exist $\ga_{1},$ $\dots,$ $\ga_{k}\in$ $ \Ga$,
$k=|G/H_{\vi_{i}}|$, such that the $k\times k$ matrix
$\big(\ga_{j}^{\vi_{i}^{g}}\big) _{\substack{j=1,\dots,k;\\
g\in G/H_{\vi_{i}}}}$ is non-degenerated. Then we set
$V_{i}=(\ga_{1},\dots,\ga_{k})$. Since the multiplication on
$\ga\in\Ga,\ga\ne0$ does not change the support, we obtain
$$|\supp y|\leq k|\supp x|\leq |G|\,|\supp x|.$$ As
$d$ we can
choose the maximum of   $d_{i}=1+\max\{ \deg v\mid v\in
V_{i}\}$, $i=1,\dots,n$. It proves
\eqref{enumerate-is-in-support}.

Now we prove \eqref{enumerate-extracted-from-support}. If
$\supp u=G\cdot m$ then $u=[bm]$ for some $b\in L^{H_{m}} $ and there is nothing to
prove. Fix some $k\geq 2$. Assume $u=[c m]+v,m\not\in\supp
v, |\supp v|\leq k-1$. For $f\in\Gamma_{1}$ consider the
polynomial $f_{S}$, (subsection
\ref{subsection-Properties-of-Galois-algebras},
\eqref{equation-relations-in-as-in-gz}) with  $S=\supp
u\setminus G\cdot m$. Applying
Lemma~\ref{lemma-relation-as-in-gz-acts-zero} we obtain the element
$$f_{S}\cdot u= f_{S}\cdot[cm]=\big[a\prod_{s\in
S}(f-f^{m s^{-1}}) m\big].$$
Since nonunit elements $ms^{-1}, s\in S$ act nontrivially on $\Ga$, there exists $f\in\Ga^{1}$ such that $f_{S}\cdot u$ is nonzero. Then
\begin{equation*}
[bm]:= f_{S}\cdot
u=\sum_{i=0}^{|S|} T_{i}u f^{|S|-i} , \text{ where }
T_{i}=\sum_{\substack{{T\subset S,}\\ {T=\{t_{1},\dots,t_{i}\},
 }}} f^{t_{1}}\dots f^{t_{i}}\in \Ga. \end{equation*}

Due to Condition \ref{condition-on-the-global-action} all
$f^{t},t\in S$ belong to a finite dimensional space $V$
generated by $\{\psi \Gamma_{1}\mid \psi\in \cM\}\subset
\bar{\Ga}$. Hence all $T_{i}$-th belong to the finite
dimensional space $V(k)=\Gamma\cap\ds\sum_{i=0}^{k}V^{i}$.
Denote  $C_{k}$ the maximal degree of elements from $V(k)$.
Then
\begin{equation*}
\deg [bm]\leq \max\{\deg T_{i}uf^{|S|-i}\mid i=0,\dots,|S|\}\leq C_{k}+j+|S|.
\end{equation*}
Hence we can set $t(k)=k +C_{k}.$

To prove \eqref{enumerate-globally-presented-in-support}
consider $x=[cm]\in U_{N(i)}$, $m\in \cM_i$.  By
\eqref{enumerate-is-in-support} for given
$\vi_{i}\in\cM_{1}$ there exists  $y\in U_{N(i)+d}$ such
that $\vi m\in \supp y$ and $\supp y \leq |G|$. Then by
\eqref{enumerate-extracted-from-support}
$U_{N(i)+d+t(|G|)}$ contains an element of the form
$[b\vi_{i} m]$.
\end{proof}

Now we are in the position to prove Theorem
\ref{theorem-GK-dimension-of-some-Galois}.
 Let $D=d+t(|G|)$. The space $U_{1}$
contains elements  $[a_{i}\vi_{i}]$, where $\vi_{i}$ runs
over $\cM_{1}/G$. Then, by Lemma
\ref{lemma-main-shift-in-growth},
\eqref{enumerate-globally-presented-in-support},
$U_{D(N-1)+1}$ contains a set of the form
$\tilde{\cM}_{N}=\{[c_{m}m]\,|\,m\in\cM_{N},c_{m}\ne0 \}$,
hence $U_{D(N-1)+N+1}$ contains $\Ga_{N}\tilde{\cM}_{N}$.
All elements from  $\Ga_{N}\tilde{\cM}_{N}$ are linearly
independent. But the set $B_{N}(\Ga[\cM]/G)$ is embedded
into $\Ga_{N}\tilde{\cM}_{N}$ by setting $\gamma[w]\mapsto
\gamma[c_{w}w]$, $\gamma\in \Ga_N$, $w\in \cM_{N+1}$.
Therefore,
$$d_U(N(D+1)+1-D)\geq |B_{N}(\Ga[\cM]/G)|\geq \dfrac{1}{|G|}|B_N(\Ga[\cM])|.$$
It remains to set $p=D+1, q=1-D$, $C=\dfrac{1}{|G|}$ and
apply
Lemma~\ref{lemma-on-an-estimate-of-growth-of-Galois-algebra}.

\section{Examples of Galois rings and orders}
\label{section-Examples-of-Galois-algebras}

\subsection{Generalized Weyl algebras}
\label{subsection-Generalized-Weyl-algebras}

Let $\sigma$ be an automorphism of $\Ga$ of infinite order,
$X$ and $Y$ generators of the bimodules $\Ga_{\sigma^{-1}}$
and $\Ga_{\sigma}$ respectively, $V=\Ga_{\sigma^{-1}}\oplus
\Ga_{\sigma}$, $G=\{e\}$ and $\cM$ is the cyclic group generated by
$\sigma$.  Consider a Galois order $U$ in
$K*\cM$ which is the image of some homomorphism $\tau:\Ga[V]\myto
K*\cM$ of the form $\tau(X)=a_{X} b_{X}^{-1}\sigma^{-1}$,
$\tau(Y)= a_{Y} b_{Y}^{-1}\sigma$ for some $a_{X},
b_{X},a_{Y}, b_{Y}\in \Ga\setminus\{0\}$. We can assume
$a_{X}= b_{X}=1$.  The element $a=a_{Y}b_{Y}^{-1}$ defines
a 2-cocycle  $\xi:\bZ\times \bZ\rightarrow K^{*}$, such
that $\xi(-1,1)=a$.
The following statement is obvious.

\begin{proposition}
\label{proposition-gen-GWA}  $U$ is a Galois
order over $\Ga$ if and only if $a\in \Ga$. In this case $U$ is isomorphic to
 a generalized Weyl algebra of rank $1$ (\cite{ba}), i.e. the algebra generated
over $\Ga$ by $X,Y$ subject to the relations
\begin{equation*}
X\lambda=\lambda^{\sigma}X, \quad \lambda Y=Y\lambda^{\sigma}, \quad \lambda\in \Lambda; \
YX=a,\quad
XY=a^{\sigma}.\end{equation*}
\end{proposition}

\subsection{Filtered algebras}\label{subsection-filtered-algebras}
 Let $U$ be an associative filtered
algebra over ${\k}$.

\begin{theorem}\label{thm-GKdim-PBW algebras}
Suppose $U$ is generated by  $u_1, \ldots, u_k$ over $\Ga$, $\gr U$  a polynomial ring in $N$ variables, $\cM\subset
\Aut L$ a group and $f:U\rightarrow \cK$ a homomorphism such
that $\underset{i}\cup\supp f(u_i)$ generates  $\cM$.
If $$\gkdim
\Ga+ \gro{\cM}=N$$ then $f$ is an embedding and $U$ is a Galois ring over $\Ga$.

\end{theorem}

\begin{proof} Note that
$f(U)$ is a Galois $\Ga$-ring by
Proposition~\ref{proposition-characterization-Galois-by-generators}. Also  $$\gkdim f(U)\geq \gkdim \Ga+
\gro{\cM}=N$$ by
Theorem~\ref{theorem-GK-dimension-of-some-Galois}.
Hence it is enough to prove that
$I=\Ker f$ equals zero. Assume $I\neq 0$. Then
$$N=\gkdim U=\gkdim \gr U>\gkdim \gr U/\gr I=\gkdim f(U)\geq N$$ which is a contradiction.
\end{proof}

Below in \ref{gl-Galois} Theorem~\ref{thm-GKdim-PBW algebras} will be applied to construct examples of Galois rings.

\subsubsection{ General linear Lie algebras}
 \label{gl-Galois}

Let $\gl_n$ be the general linear Lie algebra over $\k$,
$e_{ij}, i,j=1, \ldots, n$ its standard basis, $U_{n}$ its
universal enveloping algebra and $Z_{n}$ the center of
$U_{n}$. Then we have natural embeddings on the left upper corner
$$\gl_1\subset \gl_2\subset \ldots \subset \gl_n\text{ and induced embeddings } U_1\subset U_2\subset \ldots\subset U_n.$$
The \emph{Gelfand-Tsetlin} subalgebra
${\Gamma}$ in $U_{n}$ is generated by $\{
Z_m\,|\,m=1,\ldots, n \}$, which is a polynomial algebra in $\frac{n(n+1)}{2}$ variables.
 Denote by $K$ be
the field of fractions of ${\Gamma}$.
In the paper \cite{zh:cg} was constructed a system of generators $\{\lambda_{ij}\,|$ $1\leqslant
j\leqslant i\leqslant n \}$ of ${\Gamma}$ and the Galois extension $\Lambda\supset \Ga$ with the following properties.

\begin{enumerate}
\label{enumerate-properties-of-lambda}
\item  $\Lambda$ is the algebra of polynomial functions on
$\cL$ algebra in variables $\{\lambda_{ij}\mid \ell_{ij}\in \k,  1\leqslant
j\leqslant i\leqslant n \}$,   ${\mathscr L}=\Specm \Lambda$. An element $\ell=(\lambda_{ij}-\ell_{ij}\mid \ell_{ij}\in \k,  1\leqslant
j\leqslant i\leqslant n)$ of $\cL$ is usually written in the form of tableaux consisting of $n$ rows

\begin{align}\label{align-gt-tableau}
&\qquad\ell_{n1}\qquad\ell_{n2}
\qquad\qquad\cdots\qquad\qquad\ell_{nn}\nonumber\\
&\qquad\qquad\ell_{n-1,1}\qquad\ \ \cdots\ \
\ \ \qquad\ell_{n-1,n-1}\nonumber\\
&\quad\qquad\qquad\cdots\qquad\cdots\qquad\cdots\\
&\quad\qquad\qquad\qquad\ell_{21}\qquad\ell_{22}\nonumber\\
&\quad\qquad\qquad\qquad\qquad\ell_{11}  \nonumber
\end{align}

\item The product of the symmetric groups $G=$
$\overset{n}{\underset{i=1}\prod} S_i$ acts naturally on ${\mathscr
L}$, where every $S_{i}$ permutes elements of $i$-th row. This action induces the action of $G$ on $\La$.

\item  ${\Gamma}$ is identified with the invariants
$\Lambda^G$, such that  $\ga_{ij}=\sigma_{ij}(\ga_{i1},\dots,\ga_{ii})$ where $\sigma_{ij}$ is the $j$th symmetrical polynomial in $i$ variables. Denote by $L$ the fraction field
of $\Lambda$.  Then $L^G=K$ and $G=G(L/K)$ is the Galois group of the field extension $ K\subset L$.

\item Denote by $\delta_{ij}\in \cL$ a tableau whose $ij$-th element equals $1$ and all other elements are $0$. Let  $\cM \simeq$ $\mathbb
Z^{\frac{n(n-1)}{2}}$ be additive free abelian  group with free generators
$\delta^{ij}$, $1\leqslant j\leqslant i\leqslant n-1$. Analogously to \eqref{align-gt-tableau} the elements of $\cM$ are written as tableaux. Then $\cM$ acts on $\cL$ by shifts: $\delta^{ij}\cdot\ell=$ $\ell+\delta^{ij}$,
$\delta^{ij}\in$ $\cM$.  This action of $\cM$ on $\cL$ induces the action on $\Lambda$ and $L$,
hence we can consider $\cM$ as a subgroup in $\Aut L$. Note that
$G$ acts on $\cM$ by conjugations.
As in Section \ref{section-Galois-algebras} denote $\cK=(L*\cM)^{G}$.

\end{enumerate}

In  \cite{zh:cg}, Ch. X.70, Theorem 7,  the Gelfand-Tsetlin formulae (in Zhelobenko form) are given for the action of generators of $\gl_{n}$ on a Gelfand-Tsetlin basis of a finite dimensional irreducible representation.
We show that these formulae in fact endow $U(\gl_{n})$ with a structure of a Galois  order (Proposition \ref{proposition-galois-realization-gl}).
We need the following corollary from the Gelfand-Tsetlin formulae (see \cite{ble} or \cite{dfo:hc}).

\begin{theorem}
\label{theorem-britten-lemair-big-series}
Let $\Omega\subset \cL$ be a set of tableaux $\ell=(\ell_{ij})$ such that $\ell_{ij}-\ell_{i'j'}\not\in\bZ$ for all possible pairs $i,i',j,j'$, $(i,j)\ne (i',j')$. Consider a $\k$-vector space $T_\ell$ with the basis $\cM$ and with the action of $E_k^+=e_{k,k+1}, E_k^{-}=e_{k+1,k}$, $k=1,\ldots,n-1$, given by the formulae
\begin{equation}
\label{equation-gz-formulae-defines-generic-module}
E_k^{\pm}\cdot m =
\displaystyle {\sum_{i=1}^k}
a_{ki}^{\pm}(\ell) (m \pm
\delta^{ki}) ,$$ where $m\in\cM$ 
 and
$$
a_{ki}^{\pm}(\ell)= \mp \frac
{\prod_j(\ell_{k\pm 1,j}-
\ell_{ki})} {\prod_{j\ne
i}(\ell_{kj}-\ell_{ki})}.
\end{equation}
The
action of an element $\ga\in\Ga$ on the basis vector $\lr{\ell}$ is just the multiplication on  $\ga(\ell)\in\k$.

\end{theorem}

Analogously to  \cite{o} we show, that the formulae \eqref{align-def-t} defines a homomorphism of $U_n$ to $\cK$.

\begin{proposition}
\label{proposition-galois-realization-gl} $U_{n}$ is a Galois ring
over $\Gamma$. This structure is defined by the embedding $\imath:U\myto
\cK$ where
\begin{align} \label{align-def-t}
&\imath (e_{k\,k+1})=\sum_{i=1}^k
\delta^{ki}a_{ki}^+=\lr{\delta^{k1}a_{k1}^+},\,\,\, \imath(e_{k+1\,k})=\sum_{i=1}^k
(-\delta^{ki})a_{ki}^-=\lr{(-\delta^{k1})a_{k1}^-},\ \\ \nonumber
&a_{ki}^{\pm}= \mp \frac
{\prod_j(\lambda_{k\pm 1,j}-
\lambda_{ki})} {\prod_{j\ne
i}(\lambda_{kj}-\lambda_{ki})}, \text{ for } k=1,\dots, n.
\end{align}
\end{proposition}

\begin{proof}

Let $S$ be the multiplicative $\cM$-invariant subset in $\Ga$, generated by $\la_{ij}-\la_{ij'}-k$ for all possible $i,i',j,j'$ with $(i,j)\ne(i',j')$, where $k$ running $\bZ$, and  $\La_{S}$ the corresponding localization. Then $\La_{S}*\cM$ has a structure of a $\La_{S}*\cM$-bimodule and for every  $\ell\in\Omega=\Specm \La_{S}$ is defined a left $\La_{S}*\cM$-module $$V_{\ell}= (\La_{S}*\cM)\otimes_{\La_{S}}(\La_{S}/ \ell).$$

Analogously the action from the left by elements $\displaystyle {\sum_{i=1}^k}   (\pm\delta^{ki})
a_{k i}^{\pm}(\la), k=1,\dots,n-1$ defines on $V(\ell)$ the structure of the left $U$-module, isomorphic to the module $T_{\ell}$
from Theorem \ref{theorem-britten-lemair-big-series}. These module structures define homomorphisms of $\k$-algebras $$\tau_{\ell}:U\myto \End_{\k}(V_{\ell})\textrm{ and } \rho_{\ell}:\La_{S}*\cM\myto \End_{\k}(V_{\ell}),$$ besides $\Im \tau_{\ell}\subset
\Im \rho_{\ell}$. Hence we have the diagonal homomorphisms of $\k$-algebras
$$\Delta_{\tau}:U\myto \ds\prod_{\ell\in\Omega} \End_{\k}(V_{\ell}) \textrm{ and }
\Delta_{\rho}:\La_{S}*\cM\myto \ds\prod_{\ell\in\Omega} \End_{\k}(V_{\ell}),$$
 again $\Im \Delta_{\tau}\subset \Im \Delta_{\rho}$.
But $\Delta_{\rho}$ is an embedding, since for every nonzero $x\in \La_{S}*\cM$ there exists $V_{\ell}$, such that $x\cdot V_{\ell}\ne0$. Hence the mappings \eqref{align-def-t} from Proposition \ref{proposition-galois-realization-gl} defines the homomorphism $i:U\myto \La_{S}*\cM$. Note, that the elements in \eqref{align-def-t} belongs to $\cK$, hence $i$ defines $\imath:U\myto \cK$. To prove, that $U$ is a Galois ring note, that  $U=U(\gl_n)$ is a filtered algebra,
$\gkdim U=n^2$ and $$\gkdim \Ga + \gro{\cM}=\frac{n(n+1)}{2}+\frac{n(n-1)}{2}=n^{2}. $$ Applying
Theorem~\ref{thm-GKdim-PBW algebras} we conclude that
$\imath$ is an embedding and thus $U$ is a Galois ring.

Now we give here two different proofs of the fact that $U$ is a Galois order.

First method to prove that $U=U(\gl_{n})$ is a Galois order
is based on Proposition~\ref{prop-reduction-to-commutative}. Let
$X=(x_{ij})$ be $n\times n$-matrix with indeterminates $x_{ij}$,
$X_k$ its submatrix of size $k\times k$, formed by the intersection of  the first $k$ rows and the first $k$ columns of $X$,  $\chi_{ki}$ ($i\leq k$)
$i$-th coefficient of the characteristic polynomial of $X_k$.
In the case of $U(\gl_{n})$ corresponding graded algebra $\bar{U}$ can be identified with the  polynomial algebra in the variables $x_{ij}$, $1\leq i,j\leq n$ and the image of the canonical embedding $\imath: \gr{\Gamma} \hookrightarrow \gr{U}$ (see Proposition \ref{prop-reduction-to-commutative}) is generated by $\chi_{ki}, 1\leq k\leq n; 1\leq i\leq k $.
The $\Specm \gr{U}$ in a natural way can be interpreted as  the space $n\times n$ matrices.  Besides the induces map $\imath^{*}:\Specm \gr{U} \to
\Specm \gr{\Gamma}$
 is the map $$\bC^{n^2}\myto \bC^{n(n+1)/2},
A \longmapsto\big(\chi_{ki}(A_{k})\mid k=1, \ldots, n; i=1, \ldots, k\big),$$ defined in \cite{kw}.  It is known, that this map  is an
epimorphism (\cite{kw},  Theorem 1). Then Proposition~\ref{prop-reduction-to-commutative} implies that $U$ is a Galois order.

Another method is based on the paper \cite{o1}, where is was shown that the variety $(\imath^{*})^{-1}(0)$ is an equidimensional variety of dimension $\dfrac{n(n-1)}{2}$.
Further, from this fact in \cite{fo1} it is deduced that $U$ is  free (both right and left) $\Ga$-module.   Applying now Corollary~\ref{corollary-integral-is-projective}
we conclude that $U(\gl_{n})$ is a Galois order.
\end{proof}

Realization of $U(\gl_{n})$ as a Galois order has some interesting consequences, in particular, the decomposition $\cK\simeq \bosu{\vi\in\cM/G}{}V(\vi)$ of the localization $\cK$ of $U_{n}$ by $\Ga\setminus \{0\}$; structure of the tensor category generated by $V(\vi)$´s, etc. These results will be discussed elsewhere.

\begin{remark}\label{remark-embedding-gl-by}
Realization of $U(\gl_{n})$ as a Galois order is analogous to the
embedding of $U(\gl_{n})$ into a product of localized Weyl algebras
constructed in \cite{khm}.
\end{remark}

\begin{remark}
\label{remark-finite-W-algebras}
The developed techniques can be used effectively in the case of finite $W$-algebras.
Let $\mathfrak g=\gl_m$, $f\in \mathfrak g$, $\mathfrak g=\oplus^{}_{j\in \mathbb Z}\mathfrak g_j$
 a \emph{good grading} for $f$, i.e.
$f\in \mathfrak g_2$ and $ad \, f$ is injective
on $\mathfrak g_j$ for $j\leqslant -1$ and surjective for
$j\geqslant -1$. A non-degenerate invariant symmetric bilinear form $(.\,\, ,.)$ on $\mathfrak g$
induces a non-degenerate skew-symmetric form on
$\mathfrak g_{-1}$
defined by $\langle x,y\rangle=([x,y],\, f)$.
Let $\mathcal{I}\subset {\mathfrak g}_{-1}$ be a
maximal isotropic subspace and
set $\mathfrak{t}=\bigoplus_{j\leqslant
-2}\mathfrak g_j \oplus\mathcal{I}$. Let
$\chi:U(\mathfrak{t})\rightarrow \textbf{}C$
be the one-dimensional representation such that
$x \mapsto (x,f)$
for any $x\in \mathfrak{t}$, $I_{\chi}=\Ker \chi$ and
$Q_{\chi}=U(\mathfrak g)/U(\mathfrak g)I_{\chi}$. Then
$$
\End_{U(\mathfrak g)}(Q_{\chi})^{op}.
$$
is the finite
$W$-algebra associated to the nilpotent element
$f\in \mathfrak g$.

It was shown in \cite{bk} that any finite $W$-algebra (of type $A$) is
 isomorphic to a
certain quotient of the \emph{shifted Yangian}. It is parametrized by a sequence
$\pi=(p_1,
ldots, p_n)$ with
$p_1\leqslant \dots\leqslant p_n$.  We denote the corresponding $W$-algebra by $W(\pi)$.   Let  $\pi_k=(p_1,\dots,p_k)$,  $k\in\{1,\dots,n\}$. Then
we have the chain of  subalgebras
$$
W(\pi_1)\subset W(\pi_2) \subset\dots\subset W(\pi_n)=W(\pi).
$$
Denote by $\Gamma$ the  subalgebra
of $W(\pi)$ generated by the centers of
$W(\pi_k)$ for $k=1,\dots,n$.

\begin{theorem}[\cite{fmo},Theorem 6.6]
$W(\pi)$ is a Galois order over $\Gamma$.

\end{theorem}
\end{remark}

\subsection{Rings of invariant differential operators}\label{subsection-rings-invariant-diff-operators}
In this section we construct some Galois rings of invariant
differential operators on $n$-dimensional torus $\k^n\setminus
\{0\}$. Let $A_1$ be the first Weyl algebra over $\k$ generated by
$x$ and $\partial$ and $\tilde{A}_1$ its localization by $x$. Denote
$t=\partial x$. Then
$$\tilde{A}_1\simeq \k[t, \sigma^{\pm 1}]\simeq \k[t]*\bZ,$$
where $\sigma\in \Aut \k[t]$, $\sigma(t)=t-1$ and the first isomorphism is given by:
$x\mapsto \sigma$, $\partial \mapsto t\sigma^{-1}$. Let $\tilde{A}_n$ be the $n$-th tensor power of
$\tilde{A}_1$,
$$\tilde{A}_n\simeq \k[t_1, \ldots, t_n, \sigma_1^{\pm 1}, \ldots, \sigma_n^{\pm 1} ]\simeq \k[t_1, \ldots, t_n]*\bZ^n,$$ where $x_i, \partial_i$ are natural generators of the $n$-th Weyl algebra $A_n$, $t_i=\partial_i x_i$, $\sigma_i(t_j)=t_j-\delta_{ij}$, $i=1, \ldots, n$. Let $S=\k[t_1, \ldots, t_n]\setminus \{0\}$. Then in particular we have
$$A_n[S^{-1}]\simeq \k(t_1, \ldots, t_n)*\bZ^n.$$

\subsubsection{Symmetric differential operators on a torus}
The symmetric group $S_n$ acts naturally on $\tilde{A}_n$ by permutations.
Denote $\Ga=\k[t_1,\ldots, t_n]^{S_n}$. Then we immediately have

\begin{proposition}\label{symmetric-diff-oper-Galois}
$\tilde{A}_n^{S_n}$ is a Galois ring over $\Ga$ in $(\k(t_1,
\ldots, t_n)*\bZ^n)^{S_n}$, where $\bZ^n$ acts on the field of
rational functions by corresponding shifts.
\end{proposition}

\subsubsection{Orthogonal differential operators on a torus}
The algebra $\tilde{A}_1$ has an involution $\ve$ such that $\ve(x)=x^{-1}$ and
$\ve(\partial)=-x^2\partial$. On the other hand $\k[t]*\bZ$ has an  involution:
$\sigma\mapsto \sigma$, $t\mapsto 2-t$. Then $\tilde{A}_1$ and $\k[t]*\bZ$ are isomorphic as involutive algebras
and the isomorphism is given by: $x\mapsto \sigma$, $\partial\mapsto t\sigma^{-1}+1-\sigma^{-2}$. Similarly we have an isomorphism of involutive algebras  $\tilde{A}_n\simeq \k[t_1, \ldots, t_n, \sigma_1^{\pm 1}, \ldots, \sigma_n^{\pm 1} ]$ and $\k[t_1, \ldots, t_n]*\bZ^n.$

 Let $W_n$ be the  Weyl group of the orthogonal Lie algebra $\mathcal O_n$.
If $n=2p+1$ then the group $W_{2p+1}=S_{p}\ltimes \bZ_{2}^{p}$ acts on $\tilde{A}_p$ where
$S_p$ acts  by the permutations of the components and
  the normal subgroup $\bZ_{2}^{p}$ is generated  by the  involutions
described above.
Consider a homomorphism $\tau:\bZ_{2}^{p}\rightarrow \bZ_{2}$ such that $(g_1,\ldots, g_p)\mapsto g_1+\ldots +g_p$ and and let $N=\Ker \tau\simeq \bZ_{2}^{p-1}$.
If
 $n=2p$ then $W_{2p}\simeq S_{p}\ltimes N$ with a natural action on  $\tilde{A}_p$. These actions induce an action of $W_n$ on $\k(t_1, \ldots, t_n)*\bZ^n$ for any $n$.
 Let $\Ga=\k[t_1,\ldots, t_n]^{W_n}$.
 Then we immediately have

 \begin{proposition}\label{orthogonal-diff-oper-Galois} Algebra $\tilde{A}_n^{W_n}$ of \emph{orthogonal differential operators on a torus}
 is a Galois ring over $\Ga$ in $(\k(t_1,
\ldots, t_n)*\bZ^n)^{W_n}$, where $\bZ^n$ acts on the field of
rational functions by corresponding shifts.
\end{proposition}

\subsection{Galois orders of finite rank}
\label{subsection-Galois-orders-of-finite-rank}
The following example provides a link between the theory of Galois orders and the theory of orders in the classical sense.

Let $\Lambda$ be a commutative domain integrally closed in
its fraction field $L$, $\cG\subset \Aut L$ a finite subgroup, which splits
into a semi-direct product of its subgroups $\cG=G\ltimes \cM$. Denote
$\Ga=\Lambda^{G}$ and $K=L^{G}$. Then $\Lambda$ is just the integral
closure of $\Ga$ in $L$ and the action of $G$ on $L*\cM$ is defined.
 A Galois order  $U\subset \cK$   over $\Ga$ 
will be called 
 a \emph{Galois order of finite rank}.

\begin{proposition}
\label{proposition-structure-of-Galois-orders-of-finite-rank} Let
$U\subset \cK$ be a Galois algebra of finite rank over $\Ga$ and
$E=L^{\cG}$. Then $\cK$ is a simple central algebra over $E$ and
$\dim_{E}\cK=|\cM|^{2}$.
\end{proposition}

\begin{proof}
Theorem \ref{theorem-shat-algebras}, \eqref{enum-shat-center} gives the
statement about the center, while Corollary
\ref{corollary-ideals-in-K-localization} gives the statement
about the simplicity. From
\eqref{equation-decomp-LM-G-in-invariant},
\eqref{equation-dim-K} and subsection \ref{subsection-Separation-actions} we obtain
\begin{equation}
\label{equation-K-dimension-of-K-calligraphic}
\dim_{K}\cK=\sum_{\vi\in \cM/G}
\dim_{K}(K*\cM)^{G}_{\vi}=\sum_{\vi\in  \cM/ G}
|\cO_{\vi}|=|\cM|
\end{equation}
both as a left and as a right $K$-space structure. On other hand,
$\dim_{E}K=|\cM|$, that completes the proof.
\end{proof}

\section{Acknowledgment}
\noindent The first author is supported in part by the CNPq
grant (processo 301743/2007-0) and by the Fapesp grant
(processo 2005/60337-2). The second author is grateful to
Fapesp for the financial support (processos 2004/02850-2 and 2006/60763-4)
and to the University of S\~ao Paulo for the hospitality
during his visits.

\end{document}